\documentclass[11pt,twoside]{article}
\usepackage{amsmath, amsthm, amscd, amsfonts, amssymb, graphicx, color}

\date{}

\begin{document}

\centerline{}

\centerline {\Large{\bf Continuous frames in $n$-Hilbert spaces and their tensor products }}

\newcommand{\mvec}[1]{\mbox{\bfseries\itshape #1}}
\centerline{}
\centerline{\textbf{Prasenjit Ghosh}}
\centerline{Department of Pure Mathematics, University of Calcutta,}
\centerline{35, Ballygunge Circular Road, Kolkata, 700019, West Bengal, India}
\centerline{e-mail: prasenjitpuremath@gmail.com}
\centerline{}
\centerline{\textbf{T. K. Samanta}}
\centerline{Department of Mathematics, Uluberia College,}
\centerline{Uluberia, Howrah, 711315,  West Bengal, India}
\centerline{e-mail: mumpu$_{-}$tapas5@yahoo.co.in}

\newtheorem{Theorem}{\quad Theorem}[section]

\newtheorem{definition}[Theorem]{\quad Definition}

\newtheorem{theorem}[Theorem]{\quad Theorem}

\newtheorem{remark}[Theorem]{\quad Remark}

\newtheorem{corollary}[Theorem]{\quad Corollary}

\newtheorem{note}[Theorem]{\quad Note}

\newtheorem{lemma}[Theorem]{\quad Lemma}

\newtheorem{example}[Theorem]{\quad Example}

\newtheorem{result}[Theorem]{\quad Result}
\newtheorem{conclusion}[Theorem]{\quad Conclusion}

\newtheorem{proposition}[Theorem]{\quad Proposition}

\begin{abstract}
\textbf{\emph{We introduce the notion of continuous frame in $n$-Hilbert space which is a generalization of discrete frame in $n$-Hilbert space.\,The tensor product of Hilbert spaces is a very important topic in mathematics.\,Here we also introduce the concept of continuous frame for the tensor products of $n$-Hilbert spaces.\,Further, we study dual continuous frame and continuous Bessel multiplier in $n$-Hilbert spaces and their tensor products. }}
\end{abstract}
{\bf Keywords:}  \emph{Frame, Dual frame, Continuous frame, tensor product, $n$-Hilbert space.}\\
\\
\\

\section{Introduction}
 
\smallskip\hspace{.6 cm}The notion of frame in Hilbert space was first introduced by Duffin and Schaeffer \cite{Duffin} in connection with some fundamental problem in non-harmonic analysis.\,Thereafter, it was further developed and popularized by Daubechies et al \cite{Daubechies} in 1986.\,A discrete frame is a countable family of elements in a separable Hilbert space which allows for a stable, not necessarily unique, decomposition of an arbitrary element into an expansion of the frame element.\,Continuous frames extended the concept of discrete frames when the indices are related to some measurable space.\,Continuous frame in Hilbert space was studied by A. Rahimi et al \cite{AR}.\,M. H. faroughi and E. Osgooei \cite{MH} also studied $c$-frames and $c$-Bessel mappings.\,Continuous frame and discrete frame have been used in image processing, coding theory, wavelet analysis, signal denoising, feature extraction, robust signal processing etc.

In 1970, Diminnie et al \cite{Diminnie} introduced the concept of \,$2$-inner product space.\;A generalization of \,$2$-inner product space for \,$n \,\geq\, 2$\, was developed by A.\,Misiak \cite{Misiak} in 1989.\,There are several ways to introduced the tensor product of Hilbert spaces.\,The basic concepts of tensor product of Hilbert spaces were presented by S.\,Rabinson in \cite{S} and Folland in \cite{Folland}.\\
  
In this paper, we give the notions of continuous frames in $n$-Hilbert spaces and their tensor products.\,A characterization of continuous frame in $n$-Hilbert space with the help of its pre-frame operator is discussed.\,We will see that the image of a continuous frame under a bounded invertible operator in $n$-Hilbert space is also a continuous frame in $n$-Hilbert space.\,Continuous Bessel multipliers and dual continuous frames in $n$-Hilbert spaces and their tensor product are presented.
 
\section{Preliminaries}

\begin{theorem}\cite{O}
Let \,$H_{\,1},\, H_{\,2}$\; be two Hilbert spaces and \;$U \,:\, H_{\,1} \,\to\, H_{\,2}$\; be a bounded linear operator with closed range \;$\mathcal{R}_{\,U}$.\;Then there exists a bounded linear operator \,$U^{\dagger} \,:\, H_{\,2} \,\to\, H_{\,1}$\, such that \,$U\,U^{\dagger}\,x \,=\, x\; \;\forall\; x \,\in\, \mathcal{R}_{\,U}$.
\end{theorem}

The operator \,$U^{\dagger}$\, is called the pseudo-inverse of \,$U$.

\begin{definition}\cite{AR}
Let \,$H$\, be a complex Hilbert space and \,$(\,\Omega,\, \mu\,)$\, be a measure space with positive measure \,$\mu$.\,A mapping \,$F \,:\, \Omega \,\to\, H$\, is called a continuous frame with respect to \,$\left(\,\Omega,\, \mu\,\right)$\, if
\begin{itemize}
\item[(I)] \,$F$\, is weakly-measurable, i.\,e., for all \,$f \,\in\, H$, \,$w \,\to\, \left<\,f,\, F\,(\,w\,)\,\right>$\, is a measurable function on \,$\Omega$.
\item[(II)]there exist constants \,$0 \,<\, A \,\leq\, B \,<\, \infty$\, such that
\end{itemize}
\[A\,\left\|\,f\,\right\|^{\,2} \leq \int\limits_{\Omega}\left|\,\left<\,f,\, F\,(\,w\,)\,\right>\,\right|^{\,2}\,d\mu(\,w\,) \leq B\left\|\,f\,\right\|^{\,2}\]
for all \,$f \,\in\, H$.\,The constants \,$A$\, and \,$B$\, are called continuous frame bounds.\,If \,$A \,=\, B$, then it is called a tight continuous frame.\,If the mapping \,$F$\, satisfies only the right inequality, then it is called continuous Bessel mapping with Bessel bound \,$B$.
\end{definition}

\begin{definition}\cite{MH}
Let \,$L^{\,2}\,(\,\Omega,\,\mu\,)$\, be the class of all measurable functions \,$f \,:\, \Omega \,\to\, H$\, such that \,$\|\,f\,\|^{\,2}_{\,2} \,=\,  \int\limits_{\,\Omega}\,\left\|\,f\,(\,w\,)\,\right\|^{\,2}\,d\mu(\,w\,) \,<\, \infty$.\,It can be proved that \,$L^{\,2}\,(\,\Omega,\,\mu\,)$\, is a Hilbert space with respect to the inner product defined by
\[\left<\,f,\, g\right>_{L^{\,2}} \,=\, \int\limits_{\Omega}\,\left<\,f\,(\,w\,),\, g\,(\,w\,)\right>\,d\mu(\,w\,).\]
\end{definition}

\begin{definition}\cite{MH}
Let \,$F \,:\, \Omega \,\to\, H$\, be a Bessel mapping.\;Then the operator \,$T_{C} \,:\, L^{\,2}\left(\,\Omega,\,\mu\,\right) \,\to\, H$\, is defined by
\[\left<\,T_{C}\,(\,\varphi\,),\, h\,\right> \,=\, \int\limits_{\,\Omega}\,\varphi\,(\,w\,)\,\left<\,F\,(\,w\,),\, h\,\right>\,d\mu(\,w\,)\]where \,$\varphi \,\in\, L^{\,2}\left(\,\Omega,\,\mu\,\right)$\, and \,$h \,\in\, H$\, is well-defined, linear, bounded and its adjoint operator is given by  
\[T^{\,\ast}_{C} \,:\, H \,\to\, L^{\,2}\left(\,\Omega,\,\mu\,\right) \;,\; T^{\,\ast}_{C}\,f\,(\,w\,) \,=\, \left<\,f,\, F\,(\,w\,)\,\right>\;,\; w \,\in\, \Omega.\]
The operator \,$T_{C}$\, is called a pre-frame operator or synthesis operator and its adjoint operator is called analysis operator of \,$F$.\,The operator \,$S_{C} \,:\, H \,\to\, H$\, defined by
\[\left<\,S_{C}\,(\,f\,),\, h\,\right> \,=\, \int\limits_{\,\Omega}\,\left<\,f, F\,(\,w\,)\,\right>\left<\,F\,(\,w\,),\, h\,\right>\,d\mu(\,w\,)\]is called the frame operator of \,$F$.
\end{definition}
 
\begin{definition}\cite{Upender}\label{1.def1.01}
The tensor product of Hilbert spaces \,$H$\, and \,$K$\, is denoted by \,$H \,\otimes\, K$\, and it is defined to be an inner product space associated with the inner product
\begin{equation}\label{eq1.001}   
\left<\,f \,\otimes\, g \,,\, f^{\,\prime} \,\otimes\, g^{\,\prime}\,\right> \,=\, \left<\,f,\, f^{\,\prime}\,\right>_{\,1}\;\left<\,g,\, g^{\,\prime}\,\right>_{\,2},
\end{equation}
for all \,$f,\, f^{\,\prime} \,\in\, H\; \;\text{and}\; \;g,\, g^{\,\prime} \,\in\, K$.\,The norm on \,$H \,\otimes\, K$\, is given by 
\begin{equation}\label{eq1.0001}
\left\|\,f \,\otimes\, g\,\right\| \,=\, \|\,f\,\|_{\,1}\;\|\,g\,\|_{\,2}\; \;\forall\; f \,\in\, H\; \;\text{and}\; \,g \,\in\, K.
\end{equation}
The space \,$H \,\otimes\, K$\, is complete with respect to the above inner product.\;Therefore the space \,$H \,\otimes\, K$\, is a Hilbert space.     
\end{definition} 

For \,$Q \,\in\, \mathcal{B}\,(\,H\,)$\, and \,$T \,\in\, \mathcal{B}\,(\,K\,)$, the tensor product of operators \,$Q$\, and \,$T$\, is denoted by \,$Q \,\otimes\, T$\, and defined as 
\[\left(\,Q \,\otimes\, T\,\right)\,A \,=\, Q\,A\,T^{\,\ast}\; \;\forall\; \;A \,\in\, H \,\otimes\, K.\]
It can be easily verified that \,$Q \,\otimes\, T \,\in\, \mathcal{B}\,(\,H \,\otimes\, K\,)$\, \cite{Folland}.\\

\begin{theorem}\cite{Folland}\label{th1.1}
Suppose \,$Q,\, Q^{\prime} \,\in\, \mathcal{B}\,(\,H\,)$\, and \,$T,\, T^{\prime} \,\in\, \mathcal{B}\,(\,K\,)$, then \begin{itemize}
\item[$(i)$] \,$Q \,\otimes\, T \,\in\, \mathcal{B}\,(\,H \,\otimes\, K\,)$\, and \,$\left\|\,Q \,\otimes\, T\,\right\| \,=\, \|\,Q\,\|\; \|\,T\,\|$.
\item[$(ii)$] \,$\left(\,Q \,\otimes\, T\,\right)\,(\,f \,\otimes\, g\,) \,=\, Q\,(\,f\,) \,\otimes\, T\,(\,g\,)$\, for all \,$f \,\in\, H,\, g \,\in\, K$.
\item[$(iii)$] $\left(\,Q \,\otimes\, T\,\right)\,\left(\,Q^{\,\prime} \,\otimes\, T^{\,\prime}\,\right) \,=\, (\,Q\,Q^{\,\prime}\,) \,\otimes\, (\,T\,T^{\,\prime}\,)$. 
\item[$(iv)$] \,$Q \,\otimes\, T$\, is invertible if and only if \,$Q$\, and \,$T$\, are invertible, in which case \,$\left(\,Q \,\otimes\, T\,\right)^{\,-\, 1} \,=\, \left(\,Q^{\,-\, 1} \,\otimes\, T^{\,-\, 1}\,\right)$.
\item[$(v)$] \,$\left(\,Q \,\otimes\, T\,\right)^{\,\ast} \,=\, \left(\,Q^{\,\ast} \,\otimes\, T^{\,\ast}\,\right)$.  
\end{itemize}
\end{theorem} 
 
\begin{definition}\cite{Mashadi}
A real valued function \,$\left\|\,\cdot,\, \cdots,\, \cdot \,\right\| \,:\, H^{\,n} \,\to\, \mathbb{R}$\, satisfying the following properties:
\begin{itemize}
\item[$(i)$]\;\; $\left\|\,x_{\,1},\, x_{\,2},\, \cdots,\, x_{\,n}\,\right\| \,=\, 0$\; if and only if \,$x_{\,1},\, \cdots,\, x_{\,n}$\; are linearly dependent,
\item[$(ii)$]\;\;\; $\left\|\,x_{\,1},\, x_{\,2},\, \cdots,\, x_{\,n}\,\right\|$\, is invariant under permutations of \,$x_{1},\, \cdots, x_{n}$,
\item[$(iii)$]\;\;\; $\left\|\,\alpha\,x_{\,1},\, x_{\,2},\, \cdots,\, x_{\,n}\,\right\| \,=\, |\,\alpha\,|\,\left\|\,x_{\,1},\, x_{\,2},\, \cdots,\, x_{\,n}\,\right\|$, \,$\alpha \,\in\, \mathbb{K}$,
\item[$(iv)$]\;\; $\left\|\,x \,+\, y,\, x_{\,2},\, \cdots,\, x_{\,n}\,\right\| \,\leq\, \left\|\,x,\, x_{\,2},\, \cdots,\, x_{\,n}\,\right\| \,+\,  \left\|\,y,\, x_{\,2},\, \cdots,\, x_{\,n}\,\right\|$,
\end{itemize}
for all \,$x_{\,1},\, x_{\,2},\, \cdots,\, x_{\,n},\,x,\, y \,\in\, H$, is called \,$n$-norm on \,$H$.\,A linear space \,$H$, together with a \,$n$-norm \,$\left\|\,\cdot,\, \cdots,\, \cdot \,\right\|$, is called a linear\;$n$-normed space. 
\end{definition}

\begin{definition}\cite{Misiak}
Let \,$n \,\in\, \mathbb{N}$\; and \,$H$\, be a linear space of dimension greater than or equal to \,$n$\; over the field \,$\mathbb{K}$.\;An $n$-inner product on \,$H$\, is a map 
\[\left(\,x,\, y,\, x_{\,2},\, \cdots,\, x_{\,n}\,\right) \,\longmapsto\, \left<\,x,\, y \,|\, x_{\,2},\, \cdots,\, x_{\,n} \,\right>,\; x,\, y,\, x_{\,2},\, \cdots,\, x_{\,n} \,\in\, H\]from \,$H^{n \,+\, 1}$\, to the set \,$\mathbb{K}$\, such that for every \,$x,\, y,\, x_{\,1},\, x_{\,2},\, \cdots,\, x_{\,n} \,\in\, H$,
\begin{itemize}
\item[$(i)$]\;\; $\left<\,x_{\,1},\, x_{\,1} \,|\, x_{\,2},\, \cdots,\, x_{\,n} \,\right> \,\geq\,  0$\, and \,$\left<\,x_{\,1},\, x_{\,1} \,|\, x_{\,2},\, \cdots,\, x_{\,n} \,\right> \,=\,  0$\, if and only if \,$x_{\,1},\, x_{\,2},\, \cdots,\, x_{\,n}$\, are linearly dependent,
\item[$(ii)$]\;\; $\left<\,x,\, y \,|\, x_{\,2},\, \cdots,\, x_{\,n} \,\right> \,=\, \left<\,x,\, y \,|\, x_{\,i_{\,2}},\, \cdots,\, x_{\,i_{\,n}}\,\right> $\, for every permutations \,$\left(\, i_{\,2},\, \cdots,\, i_{\,n} \,\right)$\, of \,$\left(\, 2,\, \cdots,\, n \,\right)$,
\item[$(iii)$]\;\; $\left<\,x,\, y \,|\, x_{\,2},\, \cdots,\, x_{\,n} \,\right> \,=\, \overline{\left<\,y,\, x \,|\, x_{\,2},\, \cdots,\, x_{\,n} \,\right> }$,
\item[$(iv)$]\;\; $\left<\,\alpha\,x,\, y \,|\, x_{\,2},\, \cdots,\, x_{\,n} \,\right> \,=\, \alpha \,\left<\,x,\, y \,|\, x_{\,2},\, \cdots,\, x_{\,n} \,\right> $, for \,$\alpha \,\in\, \mathbb{K}$,
\item[$(v)$]\;\; $\left<\,x \,+\, y,\, z \,|\, x_{\,2},\, \cdots,\, x_{\,n} \,\right> \,=\, \left<\,x,\, z \,|\, x_{\,2},\, \cdots,\, x_{\,n}\,\right> \,+\,  \left<\,y,\, z \,|\, x_{\,2},\, \cdots,\, x_{\,n} \,\right>$.
\end{itemize}
A linear space \,$H$\, together with an $n$-inner product \,$\left<\,\cdot,\, \cdot \,|\, \cdot,\, \cdots,\, \cdot\,\right>$\, is called an $n$-inner product space.
\end{definition}

\begin{definition}\cite{Mashadi}
A sequence \,$\{\,x_{\,k}\,\}$\, in linear\;$n$-normed space \,$H$\, is said to be convergent to \,$x \,\in\, H$\, if 
\[\lim\limits_{k \to \infty}\,\left\|\,x_{\,k} \,-\, x,\, e_{\,2},\, \cdots,\, e_{\,n} \,\right\| \,=\, 0\]
for every \,$ e_{\,2},\, \cdots,\, e_{\,n} \,\in\, H$\, and it is called a Cauchy sequence if 
\[\lim\limits_{l,\, k \,\to\, \infty}\,\left \|\,x_{l} \,-\, x_{\,k},\, e_{\,2},\, \cdots,\, e_{\,n}\,\right\| \,=\, 0\]
for every \,$ e_{\,2},\, \cdots,\, e_{\,n} \,\in\, H$.\;The space \,$H$\, is said to be complete if every Cauchy sequence in this space is convergent in \,$H$.\;An \,$n$-inner product space is called \,$n$-Hilbert space if it is complete with respect to its induce norm.
\end{definition}

\begin{definition}\label{def0.1}\cite{Prasenjit}
Let \,$H$\, be a $n$-Hilbert space and \,$a_{\,2},\, \cdots,\, a_{\,n}$\, are fixed elements in \,$H$.\;A sequence \,$\left\{\,f_{\,i}\,\right\}^{\,\infty}_{\,i \,=\, 1}$\, in \,$H$\, is said to be a frame associated to \,$\left(\,a_{\,2},\, \cdots,\, a_{\,n}\,\right)$\, if there exists constant \,$0 \,<\, A \,\leq\, B \,<\, \infty$\,  such that
\begin{equation}\label{eqqq1}
A \, \left\|\,f,\, a_{2},\, \cdots,\, a_{n} \,\right\|^{\,2} \,\leq\, \sum\limits^{\infty}_{i \,=\, 1}\,\left|\,\left<\,f,\, f_{\,i} \,|\, a_{2},\, \cdots,\, a_{n}\,\right>\,\right|^{\,2} \,\leq\, B\, \left\|\,f,\, a_{2},\, \cdots,\, a_{n}\,\right\|^{\,2}
\end{equation}  
for all \,$f \,\in\, H$.\,The constants \,$A,\,B$\, are called frame bounds.\,If \,$\left\{\,f_{\,i}\,\right\}^{\,\infty}_{\,i \,=\, 1}$\, satisfies only the right inequality of (\ref{eqqq1}), is called a Bessel sequence associated to \,$\left(\,a_{\,2},\, \cdots,\, a_{\,n}\,\right)$\, in \,$H$\, with bound \,$B$.
\end{definition}

Let \,$a_{\,2},\, a_{\,3},\, \cdots,\, a_{\,n}$\, be the fixed elements in \,$H$\, and \,$L_{F}$\, denote the linear subspace of \,$H$\, spanned by the non-empty finite set \,$F \,=\, \left\{\,\,a_{\,2} \,,\, a_{\,3} \,,\, \cdots \,,\, a_{\,n}\,\right\}$. Then the quotient space \,$H \,/\, L_{F}$\, is a normed linear space with respect to the norm, \,$\left\|\,x \,+\, L_{F}\,\right\|_{F} \,=\, \left\|\,x \,,\, a_{\,2} \,,\,  \cdots \,,\, a_{\,n}\,\right\|$\, for every \,$x \,\in\, H$.\;Let \,$M_{F}$\, be the algebraic complement of \,$L_{F}$, then \,$H \,=\, L_{F} \,\oplus\, M_{F}$.\;Define   
\[\left<\,x \,,\, y\,\right>_{F} \,=\, \left<\,x \,,\, y \;|\; a_{\,2} \,,\,  \cdots \,,\, a_{\,n}\,\right>\; \;\text{on}\; \;H.\]
Then \,$\left<\,\cdot \,,\, \cdot\,\right>_{F}$\, is a semi-inner product on \,$H$\, and this semi-inner product induces an inner product on the quotient space \,$H \,/\, L_{F}$\; which is given by
\[\left<\,x \,+\, L_{F} \,,\, y \,+\, L_{F}\,\right>_{F} \,=\, \left<\,x \,,\, y\,\right>_{F} \,=\, \left<\,x \,,\, y \,|\, a_{\,2} \,,\,  \cdots \,,\, a_{\,n} \,\right>\;\; \;\forall \;\; x,\, y \,\in\, H.\]
By identifying \,$H \,/\, L_{F}$\; with \,$M_{F}$\; in an obvious way, we obtain an inner product on \,$M_{F}$.\;Then \,$M_{F}$\, is a normed space with respect to the norm \,$\|\,\cdot\,\|_{F}$\, defined by \,$\|\,x\,\|_{F} \;=\; \sqrt{\left<\,x \,,\, x \,\right>_{F}}\; \;\forall\, x \,\in\, M_{F}$.\;Let \,$H_{F}$\, be the completion of the inner product space \,$M_{F}$ \cite{Prasenjit}.

\begin{theorem}\label{th2}\cite{Prasenjit}
Let \,$H$\, be a n-Hilbert space.\;Then \,$\left\{\,f_{\,i}\,\right\}^{\,\infty}_{\,i \,=\, 1} \,\subseteq\, H$\; is a frame associated to \,$\left(\,a_{\,2},\, \cdots,\, a_{\,n}\,\right)$\; with bounds \,$A$\, and \,$B$\; if and only if it is a frame for the Hilbert space \,$H_{F}$\; with bounds \,$A$\, and \,$B$.
\end{theorem}

For more details on frames in \,$n$-Hilbert spaces and their tensor products one can go through the papers \cite{Prasenjit,G,PK}.

\section{Continuous frame in $n$-Hilbert space}
In this section, first we give the definition of a continuous frame in $n$-Hilbert space and then discuss some of its properties. 

\begin{definition}
Let \,$H_{1}$\, be a complex \,$n$-Hilbert space and \,$a_{\,2},\,\cdots,\,a_{\,n} \,\in\, H_{1}$\, and \,$\left(\,\Omega,\, \mu\,\right)$\, be a measure space with positive measure \,$\mu$.\,A mapping \,$\mathbb{F} \,:\, \Omega \,\to\, H_{1}$\, is called a continuous frame or \,$c$-frame associated to \,$\left(\,a_{\,2},\,\cdots,\,a_{\,n}\,\right)$\, with respect to \,$\left(\,\Omega,\, \mu\,\right)$\, if
\begin{itemize}
\item[(I)] \,$\mathbb{F}$\, is weakly-measurable, i.\,e., for all \,$f \,\in\, H_{1}$, the mapping given by \,$w \,\to\, \left<\,f,\, \mathbb{F}\,(\,w\,) \,|\, a_{\,2},\,\cdots,\,a_{\,n}\,\right>$\, is a measurable function on \,$\Omega$.
\item[(II)]there exist constants \,$0 \,<\, A \,\leq\, B \,<\, \infty$\, such that
\end{itemize}
\begin{align}
A\,\left\|\,f,\,a_{\,2},\,\cdots,\,a_{\,n}\,\right\|^{\,2} &\leq \int\limits_{\,\Omega}\left|\,\left<\,f,\, \mathbb{F}\,(\,w\,) \,|\, a_{\,2},\,\cdots,\,a_{\,n}\,\right>\,\right|^{\,2}\,d\mu(\,w\,)\nonumber\\
& \leq B\left\|\,f,\,a_{\,2},\,\cdots,\,a_{\,n}\,\right\|^{\,2}\label{en1}
\end{align}
for all \,$f \,\in\, H_{1}$.\,The constants \,$A$\, and \,$B$\, are called continuous frame bounds.\,If \,$A \,=\, B$, then it is called a tight continuous frame associated to \,$\left(\,a_{\,2},\,\cdots,\,a_{\,n}\,\right)$.\,If the mapping \,$\mathbb{F}$\, satisfies only the right inequality of (\ref{en1}), then it is called continuous Bessel mapping associated to \,$\left(\,a_{\,2},\,\cdots,\,a_{\,n}\,\right)$\, with Bessel bound \,$B$.      
\end{definition}

If \,$\mu$\, is a counting measure and \,$\mu \,=\, \mathbb{N}$, \,$\mathbb{F}$\, is called a discrete frame associated to \,$\left(\,a_{\,2},\,\cdots,\,a_{\,n}\,\right)$\, for \,$H_{1}$.

\begin{note}
Let \,$\left(\,\Omega,\, \mu\,\right)$\, be a measure space with \,$\mu$\, is \,$\sigma$-finite.\,Then the mapping \,$\mathbb{F} \,:\, \Omega \,\to\, H_{1}$\, is a continuous frame associated to \,$\left(\,a_{\,2},\,\cdots,\,a_{\,n}\,\right)$\, with bounds \,$A$\, and \,$B$\, if and only if it is a continuous frame for the Hilbert space \,$H_{F}$\, with bounds \,$A$\, and \,$B$.   
\end{note}

\begin{note}
Define the representation space \,$L_{F}^{\,2}\left(\,\Omega,\, \mu\,\right)$
\[ \,=\, \left\{\,\varphi \,:\, \Omega \,\to\, H_{F} \,|\, \,\varphi\, \,\text{is measurable and}\,\int\limits_{\,\Omega}\left\|\,\varphi\,(\,w\,),\,a_{\,2},\,\cdots,\,a_{\,n}\,\right\|^{\,2}\,d\mu(\,w\,) \,<\, \infty\,\right\}.\]
It can be easily proved that \,$L_{F}^{\,2}\left(\,\Omega,\, \mu\,\right)$\, is a Hilbert space with respect to the inner product defined by
\[\left<\,\varphi,\, \psi\,\right>_{L_{F}^{\,2}} \,=\, \int\limits_{\,\Omega}\,\left<\,\varphi\,(\,w\,),\,\psi\,(\,w\,) \,|\, a_{\,2},\,\cdots,\,a_{\,n}\,\right>\,d\mu(\,w\,)\; \;\text{for}\; \;\varphi,\, \psi \,\in\, L_{F}^{\,2}\left(\,\Omega,\, \mu\,\right).\]
\end{note}

\begin{theorem}\label{th2.1}
Let \,$\left(\,\Omega,\, \mu\,\right)$\, be a measure space and \,$\mathbb{F} \,:\, \Omega \,\to\, H_{1}$\, be a continuous Bessel mapping associated to \,$\left(\,a_{\,2},\,\cdots,\,a_{\,n}\,\right)$\, with bound \,$B$.\,Then the operator \,$T_{C} \,:\, L_{F}^{\,2}\left(\,\Omega,\, \mu\,\right) \,\to\, H_{F}$\, defined by
\[\left<\,T_{C}\,(\,\varphi\,),\, f \,|\, a_{\,2},\,\cdots,\,a_{\,n}\,\right> \,=\, \int\limits_{\,\Omega}\,\varphi\,(\,w\,)\,\left<\,\mathbb{F}\,(\,w\,),\, f \,|\, a_{\,2},\,\cdots,\,a_{\,n}\,\right>\,d\mu(\,w\,)\]where \,$\varphi \,\in\, L_{F}^{\,2}\left(\,\Omega,\, \mu\,\right)$\, and \,$f \,\in\, H_{F}$, is well-defined, bounded and linear.\,The adjoint operator \,$T^{\,\ast}_{C} \,:\, H_{F} \,\to\, L_{F}^{\,2}\left(\,\Omega,\, \mu\,\right)$\,  given by 
\[\left(\,T^{\,\ast}_{C}\,\right)(\,w\,) \,=\, \left<\,f,\, \mathbb{F}\,(\,w\,) \,|\, a_{\,2},\,\cdots,\,a_{\,n}\,\right>,\; w \,\in\, \Omega\] is also bounded and \,$\left\|\,T_{C}\,\right\| \,=\, \left\|\,T^{\,\ast}_{C}\,\right\| \,\leq\, \sqrt{B}$.
\end{theorem}

\begin{proof}
It is easy to verify that \,$T_{C}$\, is well-defined and linear.\,Since \,$\mathbb{F}$\, is a continuous Bessel mapping associated to \,$\left(\,a_{\,2},\,\cdots,\,a_{\,n}\,\right)$\, with bound \,$B$, for each \,$\varphi \,\in\, L_{F}^{\,2}\left(\,\Omega,\, \mu\,\right)$\, and \,$f \,\in\, H_{F}$, we have
\begin{align*}
&\left\|\,T_{C}\,(\,\varphi\,),\, a_{\,2},\,\cdots,\,a_{\,n}\,\right\| \,=\, \sup\limits_{\left\|\,f,\,a_{\,2},\,\cdots,\,a_{\,n}\,\right\| \,=\, 1}\left|\,\left<\,T_{C}\,(\,\varphi\,),\, f \,|\, a_{\,2},\,\cdots,\,a_{\,n}\,\right>\,\right|\\
&\leq\, \sup\limits_{\left\|\,f,\,a_{\,2},\,\cdots,\,a_{\,n}\,\right\| \,=\, 1}\left(\,\int\limits_{\,\Omega}\,\left|\,\left<\,f,\, \mathbb{F}\,(\,w\,) \,|\, a_{\,2},\,\cdots,\,a_{\,n}\,\right>\,\right|^{\,2}\,d\mu(\,w\,)\,\right)^{1 \,/\, 2}\,\times\\
&\hspace{1cm}\left(\,\int\limits_{\,\Omega}\,\left|\,\varphi\,(\,w\,)\,\right|^{2}\,d\mu(\,w\,)\,\right)^{1 \,/\, 2}\,\\
&\leq\, \sqrt{B}\,\|\,\varphi\,\|_{\,2},
\end{align*}
Hence, \,$T_{C}$\, is bounded.\,On the other hand, for each \,$\varphi \,\in\, L_{F}^{\,2}\left(\,\Omega,\, \mu\,\right)$\, and \,$f \,\in\, H_{F}$, 
\begin{align*}
&\left<\,T^{\,\ast}_{C}\,(\,f\,),\, \varphi \,|\, a_{\,2},\,\cdots,\,a_{\,n}\,\right> \,=\, \left<\,f,\, T_{C}\,(\,\varphi\,) \,|\, a_{\,2},\,\cdots,\,a_{\,n}\,\right>\\
&\hspace{1cm}=\, \int\limits_{\,\Omega}\,\overline{\varphi\,(\,w\,)}\,\left<\,f,\, \mathbb{F}\,(\,w\,) \,|\, a_{\,2},\,\cdots,\,a_{\,n}\,\right>\,d\mu(\,w\,)\\
&\hspace{1cm} \,=\, \left<\,\left<\,f,\, \mathbb{F} \,|\, a_{\,2},\,\cdots,\,a_{\,n}\,\right>,\, \varphi \,|\, a_{\,2},\,\cdots,\,a_{\,n}\,\right>. 
\end{align*}
This verify that 
\[\left(\,T^{\,\ast}_{C}\,f\,\right)(\,w\,) \,=\, \left<\,f,\, \mathbb{F}\,(\,w\,) \,|\, a_{\,2},\,\cdots,\,a_{\,n}\,\right>,\; w \,\in\, \Omega.\]Also, for each \,$f \,\in\, H_{F}$, we have
\begin{align*}
\left\|\,T^{\,\ast}_{C}\,(\,f\,),\, \varphi \,|\, a_{\,2},\,\cdots,\,a_{\,n}\,\right\|^{\,2}& \,=\, \left<\,T^{\,\ast}_{C}\,(\,f\,),\,T^{\,\ast}_{C}\,(\,f\,) \,|\, a_{\,2},\,\cdots,\,a_{\,n}\,\right>\\
&=\,\int\limits_{\,\Omega}\,\left|\,\left<\,f,\, \mathbb{F}\,(\,w\,) \,|\, a_{\,2},\,\cdots,\,a_{\,n}\,\right>\,\right|^{\,2}\,d\mu(\,w\,).
\end{align*}
This implies that
\begin{align*}
\|\,T_{C}\,\| &\,=\, \sup\limits_{\left\|\,f,\,a_{\,2},\,\cdots,\,a_{\,n}\,\right\| \,=\, 1}\left(\,\int\limits_{\,\Omega}\,\left|\,\left<\,f,\, \mathbb{F}\,(\,w\,) \,|\, a_{\,2},\,\cdots,\,a_{\,n}\,\right>\,\right|^{\,2}\,d\mu(\,w\,)\,\right)^{1 \,/\, 2}\\
& \,\leq\, \sqrt{B}. 
\end{align*} 
\end{proof}

\begin{note}
The operator \,$T_{C}$\, defined in the Theorem \ref{th2.1}, is called a pre-frame operator or synthesis operator and \,$T^{\,\ast}_{C}$\, is called an analysis operator of \,$\mathbb{F}$. 
\end{note}

\begin{definition}
The operator \,$S_{C} \,:\, H_{F} \,\to\, H_{F}$\, defined by 
\[S_{C}\,(\,f\,)\,(\,w\,) \,=\, T_{C}\,T^{\,\ast}_{C}\,(\,f\,)\,(\,w\,) \,=\, T_{C}\,\left(\,\left<\,f,\, \mathbb{F}\,(\,w\,) \,|\, a_{\,2},\,\cdots,\,a_{\,n}\,\right>\,\right)\] 
\[\,=\, \int\limits_{\,\Omega}\,\left<\,f,\, \mathbb{F}\,(\,w\,) \,|\, a_{\,2},\,\cdots,\,a_{\,n}\,\right>\,F\,(\,w\,)\,d\mu(\,w\,),\]is called continuous frame operator of \,$\mathbb{F}$.
\end{definition}

\begin{remark}\label{rr1}
Let \,$\mathbb{F} \,:\, \Omega \,\to\, H_{1}$\, be a continuous frame associated to \,$\left(\,a_{\,2},\,\cdots,\,a_{\,n}\,\right)$\, for \,$H_{1}$\, with respect to \,$\left(\,\Omega,\, \mu\,\right)$.\,For each \,$f,\, g \,\in\, H_{F}$, we have
\begin{align*}
&\left<\,S_{C}\,f,\, g \,|\, a_{\,2},\,\cdots,\,a_{\,n}\,\right>\\
& \,=\, \int\limits_{\,\Omega}\,\left<\,f,\, \mathbb{F}\,(\,w\,) \,|\, a_{\,2},\,\cdots,\,a_{\,n}\,\right>\,\left<\,\mathbb{F}\,(\,w\,),\, g \,|\, a_{\,2},\,\cdots,\,a_{\,n}\,\right>\,d\mu(\,w\,).
\end{align*}
Thus, for each \,$f \,\in\, H_{F}$, we get
\begin{align*}
&\left<\,S_{C}\,f,\, f \,|\, a_{\,2},\,\cdots,\,a_{\,n}\,\right>\\
& \,=\, \int\limits_{\,\Omega}\,\left|\,\left<\,f,\, \mathbb{F}\,(\,w\,) \,|\, a_{\,2},\,\cdots,\,a_{\,n}\,\right>\,\right|^{\,2}\,d\mu(\,w\,).
\end{align*}

Therefore, for each \,$f \,\in\, H_{F}$, from (\ref{en1}), we get
\[A\,\left<\,f,\, f \,|\, a_{\,2},\,\cdots,\,a_{\,n}\,\right> \,\leq\, \left<\,S_{C}\,f,\, f \,|\, a_{\,2},\,\cdots,\,a_{\,n}\,\right> \,\leq\, B\,\left<\,f,\, f \,|\, a_{\,2},\,\cdots,\,a_{\,n}\,\right>.\]
Hence, \,$A\,I_{F} \,\leq\, S_{C} \,\leq\, B\,I_{F}$.
\end{remark}

\begin{theorem}\label{th2.2}
Let \,$\left(\,\Omega,\, \mu\,\right)$\, be a measure space, where \,$\mu$\, is a \,$\sigma$-finite measure and let \,$\mathbb{F} \,:\, \Omega \,\to\, H_{1}$\, be a measurable function.\,If the operator \,$T_{C} \,:\, L_{F}^{\,2}\left(\,\Omega,\, \mu\,\right) \,\to\, H_{F}$\, defined by
\[\left<\,T_{C}\,(\,\varphi\,),\, f \,|\, a_{\,2},\,\cdots,\,a_{\,n}\,\right> \,=\, \int\limits_{\,\Omega}\,\varphi\,(\,w\,)\,\left<\,\mathbb{F}\,(\,w\,),\, f \,|\, a_{\,2},\,\cdots,\,a_{\,n}\,\right>\,d\mu(\,w\,)\]where \,$\varphi \,\in\, L_{F}^{\,2}\left(\,\Omega,\, \mu\,\right)$\, and \,$f \,\in\, H_{F}$, is a bounded operator, then \,$\mathbb{F}$\, is a continuous Bessel mapping associated to \,$\left(\,a_{\,2},\,\cdots,\,a_{\,n}\,\right)$.  
\end{theorem}

\begin{proof}
By the Theorem \ref{th2.1}, we have 
\[T^{\,\ast}_{C}\,(\,f\,)(\,w\,) \,=\, \left<\,f,\, \mathbb{F}\,(\,w\,) \,|\, a_{\,2},\,\cdots,\,a_{\,n}\,\right>,\; w \,\in\, \Omega.\]\,Now, for each \,$f \,\in\, H_{F}$, we have
\begin{align*}
\int\limits_{\,\Omega}\,\left|\,\left<\,f,\, \mathbb{F}\,(\,w\,) \,|\, a_{\,2},\,\cdots,\,a_{\,n}\,\right>\,\right|^{\,2}\,d\mu(\,w\,)& \,=\, \left\|\,T^{\,\ast}_{C}\,(\,f\,),\, \varphi \,|\, a_{\,2},\,\cdots,\,a_{\,n}\,\right\|^{\,2}\\
&\leq\, \left\|\,T_{C}\,\right\|^{\,2}\,\left\|\,f,\,a_{\,2},\,\cdots,\,a_{\,n}\,\right\|^{\,2}.
\end{align*}
This completes the proof.
\end{proof}

In the next theorem, we give a characterization of a continuous frame associated to \,$\left(\,a_{\,2},\,\cdots,\,a_{\,n}\,\right)$\, for \,$H_{1}$\, with respect to its pre-frame operator under some sufficient conditions.

\begin{theorem}
Let \,$\left(\,\Omega,\, \mu\,\right)$\, be a measure space, where \,$\mu$\, is a \,$\sigma$-finite measure.\,Then the mapping \,$\mathbb{F} \,:\, \Omega \,\to\, H_{1}$\, is a continuous frame associated to \,$\left(\,a_{\,2},\,\cdots,\,a_{\,n}\,\right)$\, with respect to \,$\left(\,\Omega,\, \mu\,\right)$\, if and only if the pre-frame operator \,$T_{C}$\, is bounded and onto operator. 
\end{theorem}

\begin{proof}
Let \,$\mathbb{F}$\, be a continuous frame associated to \,$\left(\,a_{\,2},\,\cdots,\,a_{\,n}\,\right)$\, for \,$H_{1}$.\,Then by Theorem \ref{th2.1}, the operator \,$T_{C}$\, is bounded and it is easy to verify that \,$T_{C}$\, is one-one, onto. \\

Conversely, let \,$T_{C}$\, be bounded and onto operator.\,Then there exists a bounded operator \,$T^{\,\dagger}_{C} \,:\, H_{F} \,\to\, L_{F}^{\,2}\left(\,\Omega,\, \mu\,\right)$\, such that \,$T_{C}\,T^{\,\dagger}_{C}\,f \,=\, f\; \;\forall\; f \,\in\, H_{F}$.\,Since \,$T_{C}$\, is bounded, by Theorem \ref{th2.2}, \,$\mathbb{F}$\, is a continuous Bessel mapping associated to \,$\left(\,a_{\,2},\,\cdots,\,a_{\,n}\,\right)$\, and 
\[\left\|\,T^{\,\ast}_{C}\,(\,f\,),\, \varphi \,|\, a_{\,2},\,\cdots,\,a_{\,n}\,\right\|^{\,2} \,=\, \int\limits_{\,\Omega}\,\left|\,\left<\,f,\, \mathbb{F}\,(\,w\,) \,|\, a_{\,2},\,\cdots,\,a_{\,n}\,\right>\,\right|^{\,2}\,d\mu(\,w\,).\]
Let \,$f \,\in\, H_{F}$, then 
\[\left\|\,f,\,a_{\,2},\,\cdots,\,a_{\,n}\,\right\|^{\,2} \,\leq\, \left\|\,T^{\,\dagger}_{C}\,\right\|^{\,2}\,\left\|\,T^{\,\ast}_{C}\,(\,f\,),\, \varphi \,|\, a_{\,2},\,\cdots,\,a_{\,n}\,\right\|^{\,2}.\]
Therefore, for each \,$f \,\in\, X_{F}$, we have 
\[\left\|\,T^{\,\dagger}_{C}\,\right\|^{\,-\, 2}\left\|\,T^{\,\ast}_{C}\,(\,f\,),\, \varphi \,|\, a_{\,2},\,\cdots,\,a_{\,n}\,\right\|^{\,2} \leq \int\limits_{\,\Omega}\left|\,\left<\,f,\, \mathbb{F}\,(\,w\,) \,|\, a_{\,2},\,\cdots,\,a_{\,n}\,\right>\,\right|^{\,2}\,d\mu(\,w\,).\]
This completes the proof.  
\end{proof}

\begin{theorem}
Let \,$\mathbb{F} \,:\, \Omega \,\to\, H_{1}$\, be a continuous frame associated to \,$\left(\,a_{\,2},\,\cdots,\,a_{\,n}\,\right)$\, with respect to \,$\left(\,\Omega,\, \mu\,\right)$\, for \,$H_{1}$\, with frame operator \,$S_{C}$\, and let \,$U \,:\, H_{F} \,\to\, H_{F}$\, be a bounded and invertible operator.\,Then \,$U\,\mathbb{F}$\, is a continuous frame associated to \,$\left(\,a_{\,2},\,\cdots,\,a_{\,n}\,\right)$\, for \,$H_{1}$\, with frame operator \,$U\,S_{C}\,U^{\,\ast}$\,
\end{theorem} 

\begin{proof}
For each \,$f \,\in\, H_{F}$, we have
\[w \,\to\, \left<\,U^{\,\ast}\,f,\, \mathbb{F}\,(\,w\,) \,|\, a_{\,2},\,\cdots,\,a_{\,n}\,\right> \,=\, \left<\,f,\, U\,\mathbb{F}\,(\,w\,) \,|\, a_{\,2},\,\cdots,\,a_{\,n}\,\right>\]
is measurable.\,Since \,$U$\, is invertible, for each \,$f \,\in\, H_{F}$, we have 
\[\left\|\,f,\,a_{\,2},\,\cdots,\,a_{\,n}\,\right\| \,\leq\, \left\|\,U^{\,-\, 1}\,\right\|\,\left\|\,U^{\,\ast}\,f,\,a_{\,2},\,\cdots,\,a_{\,n}\,\right\|.\] 
Since \,$\mathbb{F}$\, is a continuous frame associated to \,$\left(\,a_{\,2},\,\cdots,\,a_{\,n}\,\right)$\, in \,$H_{1}$, for each \,$f \,\in\, H_{F}$, we have
\begin{align*}
A\,\left\|\,U^{\,\ast}\,f,\,a_{\,2},\,\cdots,\,a_{\,n}\,\right\|^{\,2} &\,\leq\, \int\limits_{\,\Omega}\,\left|\,\left<\,U^{\,\ast}\,f,\, \mathbb{F}\,(\,w\,) \,|\, a_{\,2},\,\cdots,\,a_{\,n}\,\right>\,\right|^{\,2}\,d\mu(\,w\,)\\
& \,\leq\, B\,\left\|\,U^{\,\ast}\,f,\,a_{\,2},\,\cdots,\,a_{\,n}\,\right\|^{\,2}.
\end{align*}
Therefore, for each \,$f \,\in\, H_{F}$, we have
\begin{align*}
A\,\left\|\,U^{\,-\, 1}\,\right\|^{\,-\, 2}\,\left\|\,f,\,a_{\,2},\,\cdots,\,a_{\,n}\,\right\|^{\,2}& \,\leq\, \int\limits_{\,\Omega}\,\left|\,\left<\,f,\, U\,\mathbb{F}\,(\,w\,) \,|\, a_{\,2},\,\cdots,\,a_{\,n}\,\right>\,\right|^{\,2}\,d\mu(\,w\,)\\
& \,\leq\, B\,\|\,U\,\|^{\,2}\,\left\|\,f,\,a_{\,2},\,\cdots,\,a_{\,n}\,\right\|^{\,2}.
\end{align*} 
Thus, \,$U\,\mathbb{F}$\, is a continuous frame associated to \,$\left(\,a_{\,2},\,\cdots,\,a_{\,n}\,\right)$\, with bounds \,$A\,\left\|\,U^{\,-\, 1}\,\right\|^{\,-\, 2}$\, and \,$B\,\|\,U\,\|^{\,2}$.\\

Furthermore, for each \,$f,\, g \,\in\, X_{F}$, we have
\begin{align*}
&\int\limits_{\,\Omega}\,\left<\,f,\, U\,\mathbb{F}\,(\,w\,) \,|\, a_{\,2},\,\cdots,\,a_{\,n}\,\right>\,\left<\,U\,\mathbb{F}\,(\,w\,),\, g \,|\, a_{\,2},\,\cdots,\,a_{\,n}\,\right>\,d\mu(\,w\,)\\
&=\, \int\limits_{\,\Omega}\,\left<\,U^{\,\ast}\,f,\, \mathbb{F}\,(\,w\,) \,|\, a_{\,2},\,\cdots,\,a_{\,n}\,\right>\,\left<\,\mathbb{F}\,(\,w\,),\, U^{\,\ast}\,g \,|\, a_{\,2},\,\cdots,\,a_{\,n}\,\right>\,d\mu(\,w\,)\\
&=\, \left<\,S_{C}\,U^{\,\ast}\,f,\, U^{\,\ast}\,g \,|\, a_{\,2},\,\cdots,\,a_{\,n}\,\right> \,=\, \left<\,U\,S_{C}\,U^{\,\ast}\,f,\, g \,|\, a_{\,2},\,\cdots,\,a_{\,n}\,\right>.
\end{align*}
This shows that the corresponding frame operator is \,$U\,S_{C}\,U^{\,\ast}$.
\end{proof}

Next, we end this section by discussing the continuous Bessel multiplier in \,$H_{1}$.

\begin{definition}
Let \,$\mathbb{F}$\, and \,$\mathbb{G}$\, be continuous Bessel families associated to \,$(\,a_{\,2},\, \cdots,\, a_{\,n}\,)$\, for \,$H_{1}$\, with respect to \,$\left(\,\Omega,\, \mu\,\right)$\, having bounds \,$B_{1}$\, and \,$B_{2}$\, and \,$m \,:\, \Omega \,\to\, \mathbb{C}$\, be a measurable function.\,The operator \,$M_{m,\, \mathbb{F},\, \mathbb{G}} \,:\, H_{F} \,\to\, H_{F}$\, defined by
\begin{align*}
&\left<\,M_{m,\, \mathbb{F},\, \mathbb{G}}\,f,\, g \,|\, a_{\,2},\,\cdots,\,a_{\,n}\,\right>\\
& \,=\, \int\limits_{\,\Omega}\,m\,(\,w\,)\,\left<\,f,\, \mathbb{F}\,(\,w\,) \,|\, a_{\,2},\, \cdots,\, a_{\,n}\,\right>\,\left<\,\mathbb{G}\,(\,w\,),\, g \,|\, a_{\,2},\, \cdots,\, a_{\,n}\,\right>\,d\mu\,(\,w\,)
\end{align*} 
is called continuous Bessel multiplier associated to \,$(\,a_{\,2},\, \cdots,\, a_{\,n}\,)$\, of \,$\mathbb{F}$\, and \,$\mathbb{G}$\, with respect to \,$m$. 
\end{definition}

\begin{theorem}\label{thm4.22}
The continuous Bessel multiplier associated to \,$(\,a_{\,2},\, \cdots,\, a_{\,n}\,)$\, of \,$\mathbb{F}$\, and \,$\mathbb{G}$\, with respect to \,$m$\, is well defined and bounded. 
\end{theorem}

\begin{proof}
For any \,$f,\, g \,\in\, H_{F}$, we have
\begin{align*}
&\left|\,\left<\,M_{m,\, \mathbb{F},\, \mathbb{G}}\,f,\, g \,|\, a_{\,2},\,\cdots,\,a_{\,n}\,\right>\,\right|\\
&=\,\left|\int\limits_{\,\Omega}\,m\,(\,w\,)\left<\,f,\, \mathbb{F}\,(\,w\,) \,|\, a_{\,2},\, \cdots,\, a_{\,n}\,\right>\left<\,\mathbb{G}\,(\,w\,),\, g \,|\, a_{\,2},\, \cdots,\, a_{\,n}\,\right>d\mu\,(\,w\,)\,\right|\\
&\leq\,\|\,m\,\|_{\,\infty}\left(\,\int\limits_{\,X_{1}}\,\left|\,\left<\,f,\, \mathbb{F}\,(\,w\,) \,|\, a_{\,2},\, \cdots,\, a_{\,n}\,\right>\,\right|^{\,2}d\mu\,(\,w\,)\,\right)^{1 \,/\, 2}\,\times\\
&\hspace{1cm}\left(\,\int\limits_{\,\Omega}\,\left|\,\left<\,g,\, \mathbb{G}\,(\,w\,) \,|\, a_{\,2},\, \cdots,\, a_{\,n}\,\right>\,\right|^{\,2}d\mu\,(\,w\,)\,\right)^{1 \,/\, 2}\,\times\\
&\leq\,\|\,m\,\|_{\,\infty}\,\sqrt{B_{\,1}\,B_{\,2}}\,\left\|\,f,\,a_{\,2},\,\cdots,\,a_{\,n}\,\right\|^{\,2}\,\left\|\,g,\,a_{\,2},\,\cdots,\,a_{\,n}\,\right\|^{\,2}.
\end{align*}
This shows that \,$\left\|\,M_{m,\, \mathbb{F},\, \mathbb{G}}\,\right| \,\leq\, \|\,m\,\|_{\,\infty}\,\sqrt{B_{\,1}\,B_{\,2}}$\, and so \,$M_{m,\, \mathbb{F},\, \mathbb{G}}$\, is well-defined and bounded. 
\end{proof}

\section{Continuous frame in tensor product of $n$-Hilbert spaces}

In this section, we introduce the concept of continuous frame in tensor product of $n$-Hilbert spaces and give a characterization.\,We begin this section with the concept of tensor product of \,$n$-Hilbert spaces.\\

Let \,$H_{1}$\, and \,$H_{2}$\, be two \,$n$-Hilbert spaces associated with the \,$n$-inner products \,$\left<\,\cdot,\, \cdot \,|\, \cdot,\, \cdots,\, \cdot\,\right>_{1}$\, and \,$\left<\,\cdot,\, \cdot \,|\, \cdot,\, \cdots,\, \cdot\,\right>_{2}$, respectively.\;The tensor product of \,$H_{1}$\, and \,$H_{2}$\, is denoted by \,$H_{1} \,\otimes\, H_{2}$\, and it is defined to be an \,$n$-inner product space associated with the \,$n$-inner product given by 
\begin{align}
&\left<\,f \,\otimes\, g,\, f_{\,1} \,\otimes\, g_{\,1} \,|\, f_{\,2} \,\otimes\, g_{\,2},\, \,\cdots,\, f_{\,n} \,\otimes\, g_{\,n}\,\right>\nonumber\\
& \,=\, \left<\,f,\, f_{\,1} \,|\, f_{\,2},\, \,\cdots,\, f_{\,n}\,\right>_{1}\,\left<\,g,\, g_{\,1} \,|\, g_{\,2},\, \,\cdots,\, g_{\,n}\,\right>_{2},\label{2.eq2.31}
\end{align}
for all \,$f,\, f_{\,1},\, f_{\,2},\, \,\cdots,\, f_{\,n} \,\in\, H_{1}$\, and \,$g,\, g_{\,1},\, g_{\,2},\, \,\cdots,\, g_{\,n} \,\in\, H_{2}$.\\
The \,$n$-norm on \,$H_{1} \,\otimes\, H_{2}$\, is defined by
\begin{align}
&\left\|\,f_{\,1} \,\otimes\, g_{\,1},\, f_{\,2} \,\otimes\, g_{\,2},\, \,\cdots,\,\, f_{\,n} \,\otimes\, g_{\,n}\,\right\|\nonumber\\
& =\,\left\|\,f_{\,1},\, f_{\,2},\, \cdots,\, f_{\,n}\,\right\|_{1}\;\left\|\,g_{\,1},\, g_{\,2},\, \cdots,\, g_{\,n}\,\right\|_{2},\label{2.eq2.32}
\end{align} 
for all \,$f_{\,1},\, f_{\,2},\, \,\cdots,\, f_{\,n} \,\in\, H_{1}\, \;\text{and}\; \,g_{\,1},\, g_{\,2},\, \,\cdots,\, g_{\,n} \,\in\, H_{2}$, where the \,$n$-norms \,$\left\|\,\cdot,\, \cdots,\, \cdot \,\right\|_{1}$\, and \,$\left\|\,\cdot,\, \cdots,\, \cdot \,\right\|_{2}$\, are generated by \,$\left<\,\cdot,\, \cdot \,|\, \cdot,\, \cdots,\, \cdot\,\right>_{1}$\, and \,$\left<\,\cdot,\, \cdot \,|\, \cdot,\, \cdots,\, \cdot\,\right>_{2}$, respectively.\;The space \,$H_{1} \,\otimes\, H_{2}$\, is completion with respect to the above \,$n$-inner product.\;Therefore the space \,$H_{1} \,\otimes\, H_{2}$\, is an \,$n$-Hilbert space.

Consider \,$G \,=\, \left\{\,b_{\,2},\, b_{\,3},\, \cdots,\, b_{\,n}\,\right\}$, where \,$b_{\,2},\, b_{\,3},\, \cdots,\, b_{\,n}$\, are fixed elements in \,$H_{2}$\, and \,$L_{G}$\, denote the linear subspace of \,$H_{2}$\, spanned by \,$G$.\,Now, we can define the Hilbert space \,$H_{G}$\, with respect to the inner product is given by
\[\left<\,f \,+\, L_{G}\,,\, g \,+\, L_{G}\,\right>_{G} \,=\, \left<\,f \,,\, g\,\right>_{G} \,=\, \left<\,f \,,\, g \,|\, b_{\,2} \,,\,  \cdots \,,\, b_{\,n} \,\right>_{2}; \;\forall \;\; f,\, g \,\in\, H_{2}.\]

\begin{remark}
According to the definition \ref{1.def1.01}, \,$H_{F} \,\otimes\, H_{G}$\, is the Hilbert space with respect to the inner product:
\[\left<\,f \,\otimes\, g \,,\, f^{\,\prime} \,\otimes\, g^{\,\prime}\,\right> \,=\, \left<\,f \,,\, f^{\,\prime}\,\right>_{F}\;\left<\,g \,,\, g^{\,\prime}\,\right>_{G},\]
for all \,$f,\, f^{\,\prime} \,\in\, H_{F}\; \;\text{and}\; \;g,\, g^{\,\prime} \,\in\, H_{G}$.    
\end{remark}

\begin{definition}
Let \,$(\,X,\, \mu\,) \,=\, \left(\,X_{1} \,\times\, X_{2},\, \mu_{\,1} \,\otimes\, \mu_{\,2}\,\right)$\, be the product of measure spaces with \,$\sigma$-finite positive measures \,$\mu_{\,1},\, \mu_{\,2}$\, and \,$a_{\,2} \,\otimes\, b_{\,2},\, \cdots,\, a_{\,n} \,\otimes\, b_{\,n}$\, be fixed elements in \,$H_{1} \,\otimes\, H_{2}$.\,The mapping \,$\mathcal{F} \,:\, X \,\to\, H_{1} \,\otimes\, H_{2}$\, is called a continuous frame associated to \,$\left(\,a_{\,2} \,\otimes\, b_{\,2},\, \cdots,\, a_{\,n} \,\otimes\, b_{\,n}\,\right)$\, for \,$H_{1} \,\otimes\, H_{2}$\, with respect to \,$(\,X,\, \mu\,)$\, if
\begin{itemize}
\item[$(i)$]$\mathcal{F}$\, is weakly-measurable, i.\,e., for all \,$f \,\otimes\, g \,\in\, H_{1} \,\otimes\, H_{2}$, \,$x \,=\, \left(\,x_{\,1},\, x_{\,2}\,\right) \,\to\, \left<\,f \,\otimes\, g,\, \mathcal{F}\,(\,x\,) \,|\, a_{\,2} \,\otimes\, b_{\,2},\, \cdots,\, a_{\,n} \,\otimes\, b_{\,n}\,\right>$\, is a measurable function on \,$X$.
\item[$(ii)$]there exist constants \,$A,\,B \,>\, 0$\, such that
\begin{align}
&A\,\left\|\,f \,\otimes\, g,\, a_{\,2} \,\otimes\, b_{\,2},\, \cdots,\, a_{\,n} \,\otimes\, b_{\,n}\,\right\|^{\,2}\nonumber\\
& \,\leq\, \int\limits_{\,X}\,\left|\,\left<\,f \,\otimes\, g,\, \mathcal{F}\,(\,x\,) \,|\, a_{\,2} \,\otimes\, b_{\,2},\, \cdots,\, a_{\,n} \,\otimes\, b_{\,n}\,\right>\,\right|^{\,2}\,d\mu\,(\,x\,)\nonumber\\
&\leq\,B\,\left\|\,f \,\otimes\, g,\, a_{\,2} \,\otimes\, b_{\,2},\, \cdots,\, a_{\,n} \,\otimes\, b_{\,n}\,\right\|^{\,2},\label{4.eq4.11} 
\end{align}
for all \,$f \,\otimes\, g \,\in\, H_{1} \,\otimes\, H_{2}$.\,The constants \,$A$\, and \,$B$\, are called continuous frame bounds.\,If \,$A \,=\, B$, then it is called a tight continuous frame associated to \,$\left(\,a_{\,2} \,\otimes\, b_{\,2},\, \cdots,\, a_{\,n} \,\otimes\, b_{\,n}\,\right)$.\,If the mapping \,$\mathcal{F}$\, satisfies only the right inequality of (\ref{4.eq4.11}), then it is called Bessel mapping or \,$c$-Bessel mapping associated to \,$\left(\,a_{\,2} \,\otimes\, b_{\,2},\, \cdots,\, a_{\,n} \,\otimes\, b_{\,n}\,\right)$\, with Bessel bound \,$B$.  
\end{itemize}   
\end{definition}

In the following theorem, we show that the continuous frame in $n$-Hilbert space is preserved by the tensor product.

\begin{theorem}\label{th4.11}
The mapping \,$\mathcal{F} \,=\, F_{1} \,\otimes\, F_{2} \,:\, X \,\to\, H_{1} \,\otimes\, H_{2}$\, is a continuous frame associated to \,$\left(\,a_{\,2} \,\otimes\, b_{\,2},\, \cdots,\, a_{\,n} \,\otimes\, b_{\,n}\,\right)$\, for \,$H_{1} \,\otimes\, H_{2}$\, with respect to \,$(\,X,\, \mu\,)$\, if and only if \,$F_{1}$\, is a continuous frame associated to \,$\left(\,a_{\,2},\, \cdots,\, a_{\,n}\,\right)$\, for \,$H_{1}$\, with respect to \,$\left(\,X_{1},\, \mu_{\,1}\,\right)$\, and \,$F_{2}$\, is a continuous frame associated to \,$\left(\,b_{\,2},\, \cdots,\, b_{\,n}\,\right)$\, for \,$H_{2}$\, with respect to \,$\left(\,X_{2},\, \mu_{\,2}\,\right)$\,  
\end{theorem}
 
\begin{proof}
Suppose that \,$\mathcal{F} \,=\, F_{1} \,\otimes\, F_{2}$\, is a continuous frame associated to \,$(\,a_{\,2} \,\otimes\, b_{\,2},\, \cdots,\, a_{\,n} \,\otimes\, b_{\,n}\,)$\, for \,$H_{1} \,\otimes\, H_{2}$\, with respect to \,$(\,X,\, \mu\,)$.\,Let \,$f \,\in\, H_{1} \,/\, \{\,0\,\}$\, and fix \,$g \,\in\, H_{2} \,/\, \{\,0\,\}$.\,Then \,$f \,\otimes\, g \,\in\, H_{1} \,\otimes\, H_{2}$\, and by Fubini's theorem we have
\begin{align*}
&\int\limits_{\,X}\,\left|\,\left<\,f \,\otimes\, g,\, F_{1}\,(\,x_{\,1}\,) \,\otimes\, F_{2}\,(\,x_{\,2}\,) \,|\, a_{\,2} \,\otimes\, b_{\,2},\, \cdots,\, a_{\,n} \,\otimes\, b_{\,n}\,\right>\,\right|^{\,2}\,d\mu\,(\,x\,)\\
&=\,\int\limits_{\,X_{1}}\,\left|\,\left<\,f,\, F_{1}\,(\,x_{\,1}\,) \,|\, a_{\,2},\, \cdots,\, a_{\,n}\,\right>_{1}\,\right|^{\,2}\,d\mu_{\,1}\,(\,x_{1}\,)\,\times\\
&\hspace{1cm}\int\limits_{\,X_{2}}\,\left|\,\left<\,g,\, F_{2}\,(\,x_{\,2}\,) \,|\, b_{\,2},\, \cdots,\, b_{\,n}\,\right>_{2}\,\right|^{\,2}\,d\mu_{\,2}\,(\,x_{2}\,).
\end{align*}
Therefore, for each \,$f \,\otimes\, g \,\in\, H_{1} \,\otimes\, H_{2}$, (\ref{4.eq4.11}) can be written as
\begin{align*}
&A\,\left\|\,f,\, a_{\,2},\, \cdots,\, a_{\,n}\,\right\|^{\,2}_{1}\,\left\|\,g,\, b_{\,2},\, \cdots,\, b_{\,n}\,\right\|^{\,2}_{2}\\
&\leq\,\int\limits_{\,X_{1}}\,\left|\,\left<\,f,\, F_{1}\,(\,x_{\,1}\,) \,|\, a_{\,2},\, \cdots,\, a_{\,n}\,\right>_{1}\,\right|^{\,2}\,d\mu_{\,1}\,(\,x_{1}\,)\,\times\\
&\hspace{1cm}\int\limits_{\,X_{2}}\,\left|\,\left<\,g,\, F_{2}\,(\,x_{\,2}\,) \,|\, b_{\,2},\, \cdots,\, b_{\,n}\,\right>_{2}\,\right|^{\,2}\,d\mu_{\,2}\,(\,x_{2}\,)\\
&\,\leq\,B\,\left\|\,f,\, a_{\,2},\, \cdots,\, a_{\,n}\,\right\|^{\,2}_{1}\,\left\|\,g,\, b_{\,2},\, \cdots,\, b_{\,n}\,\right\|^{\,2}_{2}.
\end{align*} 
Here we may assume that every \,$F_{1}\left(\,x_{\,1}\,\right)$\, and \,$a_{\,2},\, \cdots,\, a_{\,n}$\, are linearly independent and also every \,$F_{2}\left(\,x_{\,2}\,\right)$\, and \,$b_{\,2},\, \cdots,\, b_{\,n}$\, are linearly independent.\,Hence 
\[\int\limits_{\,X_{1}}\,\left|\,\left<\,f,\, F_{1}\,(\,x_{\,1}\,) \,|\, a_{\,2},\, \cdots,\, a_{\,n}\,\right>_{1}\,\right|^{\,2}\,d\mu_{\,1}\,(\,x_{1}\,),\]
\[ \,\int\limits_{\,X_{2}}\,\left|\,\left<\,g,\, F_{2}\,(\,x_{\,2}\,) \,|\, b_{\,2},\, \cdots,\, b_{\,n}\,\right>_{2}\,\right|^{\,2}\,d\mu_{\,2}\,(\,x_{2}\,)\]are non-zero.\,Thus from the above inequality we can write
\begin{align*}
&\dfrac{A\,\left\|\,g,\, b_{\,2},\, \cdots,\, b_{\,n}\,\right\|^{\,2}_{2}}{\int\limits_{\,X_{2}}\,\left|\,\left<\,g,\, F_{2}\,(\,x_{\,2}\,) \,|\, b_{\,2},\, \cdots,\, b_{\,n}\,\right>_{2}\,\right|^{\,2}\,d\mu_{\,2}\,(\,x_{2}\,)}\,\left\|\,f,\, a_{\,2},\, \cdots,\, a_{\,n}\,\right\|^{\,2}_{1}\\
&\leq\,\int\limits_{\,X_{1}}\,\left|\,\left<\,f,\, F_{1}\,(\,x_{\,1}\,) \,|\, a_{\,2},\, \cdots,\, a_{\,n}\,\right>_{1}\,\right|^{\,2}\,d\mu_{\,1}\,(\,x_{1}\,)\\
&\leq\,\dfrac{B\,\left\|\,g,\, b_{\,2},\, \cdots,\, b_{\,n}\,\right\|^{\,2}_{2}}{\int\limits_{\,X_{2}}\,\left|\,\left<\,g,\, F_{2}\,(\,x_{\,2}\,) \,|\, b_{\,2},\, \cdots,\, b_{\,n}\,\right>_{2}\,\right|^{\,2}\,d\mu_{\,2}\,(\,x_{2}\,)}\,\left\|\,f,\, a_{\,2},\, \cdots,\, a_{\,n}\,\right\|^{\,2}_{1}.
\end{align*}
Thus, for each \,$f \,\in\, H_{1} \,/\, \{\,0\,\}$, we have
\begin{align*}
A_{\,1}\,\left\|\,f,\, a_{\,2},\, \cdots,\, a_{\,n}\,\right\|^{\,2}_{1} &\,\leq\, \int\limits_{\,X_{1}}\,\left|\,\left<\,f,\, F_{1}\,(\,x_{\,1}\,) \,|\, a_{\,2},\, \cdots,\, a_{\,n}\,\right>_{1}\,\right|^{\,2}\,d\mu_{\,1}\,(\,x_{1}\,)\\
& \leq\,B_{\,1}\,\left\|\,f,\, a_{\,2},\, \cdots,\, a_{\,n}\,\right\|^{\,2}_{1},
\end{align*}
where 
\[A_{\,1} \,=\, \inf_{g \,\in\, H_{\,2}}\,\left\{\dfrac{A\,\left\|\,g,\, b_{\,2},\, \cdots,\, b_{\,n}\,\right\|^{\,2}_{2}}{\int\limits_{\,X_{2}}\,\left|\,\left<\,g,\, F_{2}\,(\,x_{\,2}\,) \,|\, b_{\,2},\, \cdots,\, b_{\,n}\,\right>_{2}\,\right|^{\,2}\,d\mu_{\,2}\,(\,x_{2}\,)}\,\right\},\]and
\[B_{\,1} \,=\, \sup_{g \,\in\, H_{\,2}}\,\left\{\dfrac{B\,\left\|\,g,\, b_{\,2},\, \cdots,\, b_{\,n}\,\right\|^{\,2}_{2}}{\int\limits_{\,X_{2}}\,\left|\,\left<\,g,\, F_{2}\,(\,x_{\,2}\,) \,|\, b_{\,2},\, \cdots,\, b_{\,n}\,\right>_{2}\,\right|^{\,2}\,d\mu_{\,2}\,(\,x_{2}\,)}\,\right\}.\] 
This shows that \,$F_{1}$\, is a continuous frame associated to \,$\left(\,a_{\,2},\, \cdots,\, a_{\,n}\,\right)$\, for \,$H_{1}$\, with respect to \,$\left(\,X_{1},\, \mu_{\,1}\,\right)$.\,Similarly, it can be shown that \,$F_{2}$\, is a continuous frame associated to \,$\left(\,b_{\,2},\, \cdots,\, b_{\,n}\,\right)$\, for \,$H_{2}$\, with respect to \,$\left(\,X_{2},\, \mu_{\,2}\,\right)$.\\\\
Conversely, suppose that \,$F_{1}$\, is a continuous frame associated to \,$\left(\,a_{\,2},\, \cdots,\, a_{\,n}\,\right)$\, for \,$H_{1}$\, with respect to \,$\left(\,X_{1},\, \mu_{\,1}\,\right)$\, having bounds \,$A,\, B$\, and \,$F_{2}$\, is a continuous frame associated to \,$\left(\,b_{\,2},\, \cdots,\, b_{\,n}\,\right)$\, for \,$H_{2}$\, with respect to \,$\left(\,X_{2},\, \mu_{\,2}\,\right)$\, having bounds \,$C,\,D$.\,By the assumption it is easy to very that \,$F \,=\, F_{1} \,\otimes\, F_{2}$\, is weakly measurable on \,$H_{1} \,\otimes\, H_{2}$\, with respect to \,$(\,X,\, \mu\,)$.\,Now, for each \,$f \,\in\, H_{1} \,/\, \{\,0\,\},\, \,g \,\in\, H_{2} \,/\, \{\,0\,\}$, we have
\begin{align*}
A\,\left\|\,f,\, a_{\,2},\, \cdots,\, a_{\,n}\,\right\|^{\,2}_{1} &\,\leq\, \int\limits_{\,X_{1}}\,\left|\,\left<\,f,\, F_{1}\,(\,x_{\,1}\,) \,|\, a_{\,2},\, \cdots,\, a_{\,n}\,\right>_{1}\,\right|^{\,2}\,d\mu_{\,1}\,(\,x_{1}\,)\\
& \leq\,B\,\left\|\,f,\, a_{\,2},\, \cdots,\, a_{\,n}\,\right\|^{\,2}_{1},
\end{align*} 
\begin{align*}
C\,\left\|\,g,\, b_{\,2},\, \cdots,\, b_{\,n}\,\right\|^{\,2}_{2} &\,\leq\, \int\limits_{\,X_{2}}\,\left|\,\left<\,g,\, F_{2}\,(\,x_{\,2}\,) \,|\, b_{\,2},\, \cdots,\, b_{\,n}\,\right>_{2}\,\right|^{\,2}\,d\mu_{\,2}\,(\,x_{2}\,)\\
& \leq\,D\,\left\|\,g,\, b_{\,2},\, \cdots,\, b_{\,n}\,\right\|^{\,2}_{2}.
\end{align*}
Multiplying the above two inequalities and using Fubini's theorem we get
\begin{align*}
&A\,C\,\left\|\,f \,\otimes\, g,\, a_{\,2} \,\otimes\, b_{\,2},\, \cdots,\, a_{\,n} \,\otimes\, b_{\,n}\,\right\|^{\,2}\\
& \,\leq\, \int\limits_{\,X}\,\left|\,\left<\,f \,\otimes\, g,\, \mathcal{F}\,(\,x\,) \,|\, a_{\,2} \,\otimes\, b_{\,2},\, \cdots,\, a_{\,n} \,\otimes\, b_{\,n}\,\right>\,\right|^{\,2}\,d\mu\,(\,x\,)\\
&\leq\,B\,D\,\left\|\,f \,\otimes\, g,\, a_{\,2} \,\otimes\, b_{\,2},\, \cdots,\, a_{\,n} \,\otimes\, b_{\,n}\,\right\|^{\,2}, 
\end{align*} 
for all \,$f \,\otimes\, g \,\in\, H_{1} \,\otimes\, H_{2}$.\,This completes the proof.  
\end{proof} 

\begin{note}
Let \,$(\,X,\, \mu\,) \,=\, \left(\,X_{1} \,\times\, X_{2},\, \mu_{\,1} \,\otimes\, \mu_{\,2}\,\right)$\, be the product of measure spaces with \,$\sigma$-finite positive measures \,$\mu_{\,1},\, \mu_{\,2}$.\\

Let \,$L_{F \,\otimes\, G}^{\,2}\left(\,X,\, \mu\,\right)$\, be the class of all measurable functions \,$\Psi \,:\, X \,\to\, H_{F} \,\otimes\, H_{G}$\, such that 
\[\int\limits_{\,X}\left\|\,\Psi\,(\,x\,),\,a_{\,2} \,\otimes\, b_{\,2},\, \cdots,\, a_{\,n} \,\otimes\, b_{\,n}\,\right\|^{\,2}\,d\mu(\,x\,) \,<\, \infty,\]
with the inner product 
\begin{align*}
\left<\,\Psi,\, \Phi\,\right>_{L_{F \,\otimes\, G}^{\,2}} &\,=\, \int\limits_{\,X}\,\left<\,\Psi\,(\,x\,),\,\Phi\,(\,x\,) \,|\, a_{\,2} \,\otimes\, b_{\,2},\, \cdots,\, a_{\,n} \,\otimes\, b_{\,n}\,\right>\,d\mu(\,x\,),\\
&=\,\int\limits_{\,X_{1}}\,\left<\,\varphi_{1}\,(\,x_{1}\,),\,\psi_{1}\,(\,x_{1}\,) \,|\, a_{\,2},\,\cdots,\,a_{\,n}\,\right>\,d\mu(\,x_{1}\,)\,\times\\
&\hspace{1cm}\int\limits_{\,X_{2}}\,\left<\,\varphi_{2}\,(\,x_{2}\,),\,\psi_{2}\,(\,x_{2}\,) \,|\, b_{\,2},\,\cdots,\,b_{\,n}\,\right>\,d\mu(\,x_{2}\,)\\
&=\,\left<\,\varphi_{1},\, \psi_{1}\,\right>_{L_{F}^{\,2}}\,\left<\,\varphi_{2},\, \psi_{2}\,\right>_{L_{G}^{\,2}},
\end{align*}
for \,$\Psi \,=\, \varphi_{1} \,\otimes\, \varphi_{2},\, \Phi \,=\, \psi_{1} \,\otimes\, \psi_{2} \,\in\, L_{F \,\otimes\, G}^{\,2}\left(\,X,\, \mu\,\right)$.\,The space \,$L_{F \,\otimes\, G}^{\,2}\left(\,X,\, \mu\,\right)$\, is completion with respect to the above inner product.\,Therefore it is an Hilbert space.

\end{note}

\begin{remark}
Let \,$\mathcal{F}$\, be a continuous Bessel family associated to \,$(\,a_{\,2} \,\otimes\, b_{\,2},\, \cdots,\, \\a_{\,n} \,\otimes\, b_{\,n}\,)$\, for \,$H_{1} \,\otimes\, H_{2}$\, with respect to \,$(\,X,\, \mu\,)$.\,Then the synthesis operator  
\,$T_{\mathcal{F}} \,:\, L_{F \,\otimes\, G}^{\,2}\left(\,X,\, \mu\,\right) \,\to\, H_{F} \,\otimes\, H_{G}$\, defined by
\begin{align*}
T_{\mathcal{F}}\,(\,\varphi\,) &\,=\, \int\limits_{\,X}\,\varphi\,(\,x\,)\,F\,(\,x\,)\,d\mu(\,x\,)\\
& \,=\, \int\limits_{\,X_{1}}\,\int\limits_{\,X_{2}}\,\varphi\left(\,x_{\,1},\, x_{\,2}\,\right)\,F\,\left(\,x_{\,1},\, x_{\,2}\,\right)\,d\mu\left(\,x_{\,1},\, x_{\,2}\,\right)
\end{align*}
where \,$\varphi \,\in\, L_{F \,\otimes\, G}^{\,2}\left(\,X,\, \mu\,\right)$\, is well-defined, bounded and linear.\,The analysis operator \,$T^{\,\ast}_{\mathcal{F}} \,:\, H_{F} \,\otimes\, H_{G} \,\to\, L_{F \,\otimes\, G}^{\,2}\left(\,X,\, \mu\,\right)$\,  given by 
\[\left(\,T^{\,\ast}_{\mathcal{F}}\,(\,f \,\otimes\, g\,)\,\right)(\,x\,) \,=\, \left<\,f \,\otimes\, g,\, \mathcal{F}\,(\,x\,) \,|\, a_{\,2} \,\otimes\, b_{\,2},\,\cdots,\,a_{\,n} \,\otimes\, b_{\,n}\,\right>,\]$\; x \,\in\, X,\, f \,\otimes\, g \,\in\, H_{F} \,\otimes\, H_{G}$.\,The frame operator \,$S_{\mathcal{F}} \,:\, H_{F} \,\otimes\, H_{G} \,\to\, H_{F} \,\otimes\, H_{G}$\, is given by
\[S_{\mathcal{F}}\,(\,f \,\otimes\, g\,) \,=\, \int\limits_{\,X}\,\left<\,f \,\otimes\, g,\, \mathcal{F}\,(\,x\,) \,|\, a_{\,2} \,\otimes\, b_{\,2},\, \cdots,\, a_{\,n} \,\otimes\, b_{\,n}\,\right>\,\mathcal{F}\,(\,x\,)\,d\mu\,(\,x\,)\]
\end{remark}

The next theorem demonstrates that the continuous frame operator associated with the tensor product of two continuous frames in \,$n$-Hilbert spaces is exactly the tensor product of their respective continuous frame operators.

\begin{theorem}
Let \,$\mathcal{F} \,=\, F_{1} \,\otimes\, F_{2} \,:\, X \,\to\, H_{1} \,\otimes\, H_{2}$\, be a continuous frame associated to \,$\left(\,a_{\,2} \,\otimes\, b_{\,2},\, \cdots,\, a_{\,n} \,\otimes\, b_{\,n}\,\right)$\, for \,$H_{1} \,\otimes\, H_{2}$\, with respect to \,$(\,X,\, \mu\,)$.\,Then \,$S_{\mathcal{F}} \,=\, S_{F_{\,1}}  \,\otimes\, S_{F_{\,2}}$.  
\end{theorem}

\begin{proof}
Suppose that \,$\mathcal{F} \,=\, F_{1} \,\otimes\, F_{2}$\, is a continuous frame associated to \,$(\,a_{\,2} \,\otimes\, b_{\,2},\, \cdots,\, a_{\,n} \,\otimes\, b_{\,n}\,)$\, for \,$H_{1} \,\otimes\, H_{2}$\, with respect to \,$(\,X,\, \mu\,)$.\,Then for each \,$f \,\otimes\, g \,\in\, H_{F} \,\otimes\, H_{G}$, we have
\begin{align*}
&S_{\mathcal{F}}\,(\,f \,\otimes\, g\,)\\
 &\,=\, \int\limits_{\,X}\,\left<\,f \,\otimes\, g,\, F_{1}\,(\,x_{\,1}\,) \,\otimes\, F_{2}\,(\,x_{\,2}\,) \,|\, a_{\,2} \,\otimes\, b_{\,2},\, \cdots,\, a_{\,n} \,\otimes\, b_{\,n}\,\right>\,F_{1}\,(\,x_{\,1}\,) \,\otimes\, F_{2}\,(\,x_{\,2}\,)d\mu\,(\,x\,)\\
&=\,\int\limits_{\,X_{1}}\,\left<\,f,\, F_{1}\,(\,x_{\,1}\,) \,|\, a_{\,2},\, \cdots,\, a_{\,n}\,\right>_{1}\,F_{1}\,(\,x_{\,1}\,)\,d\mu_{\,1}\,(\,x_{1}\,)\,\,\otimes\\
&\hspace{1cm}\int\limits_{\,X_{2}}\,\left<\,g,\, F_{2}\,(\,x_{\,2}\,) \,|\, b_{\,2},\, \cdots,\, b_{\,n}\,\right>_{2}\,F_{2}\,(\,x_{\,2}\,)\,d\mu_{\,2}\,(\,x_{2}\,)\\
&=\, S_{F_{\,1}}\,f  \,\otimes\, S_{F_{\,2}}\,g \,=\, \left(\,S_{F_{\,1}}  \,\otimes\, S_{F_{\,2}}\,\right)\,(\,f \,\otimes\, g\,).
\end{align*} 
\end{proof}

\begin{theorem}
Let \,$F_{1}$\, be a continuous frame associated to \,$\left(\,a_{\,2},\, \cdots,\, a_{\,n}\,\right)$\, for \,$H_{1}$\, with respect to \,$\left(\,X_{1},\, \mu_{1}\,\right)$\, having bounds \,$A,\, B$\, and \,$F_{2}$\, be a continuous frame associated to \,$\left(\,b_{\,2},\, \cdots,\, b_{\,n}\,\right)$\, for \,$H_{2}$\, with respect to \,$\left(\,X_{2},\, \mu_{2}\,\right)$\, having bounds \,$C,\,D$.\,Then \,$A\,C\,I_{F \,\otimes\, G} \,\leq\, S_{F_{1} \,\otimes\, F_{2}} \,\leq\, B\,D\,I_{F \,\otimes\, G}$, where \,$\,I_{F \,\otimes\, G}$\, is the identity operator on \,$H_{F} \,\otimes\, H_{G}$\, and \,$S_{F_{1}},\, \,S_{F_{2}}$\, are continuous frame operators of \,$F_{1},\, \,F_{2}$, respectively.
\end{theorem}

\begin{proof}
Since \,$S_{F_{1}}$\, and \,$S_{F_{2}}$\, are continuous frame operators, we have 
\[A\,I_{F} \,\leq\, S_{F_{1}} \,\leq\, B\,I_{F},\; \;C\,I_{G} \,\leq\, S_{F_{2}} \,\leq\, D\,I_{G},\]
where \,$I_{F}$\,and \,$I_{G}$\, are the identity operators on \,$H_{F}$\, and \,$K_{G}$, respectively.\,Taking tensor product on the above two inequalities, we get
\begin{align*}
&A\,C \left(\,I_{F} \,\otimes\, I_{G}\,\right) \,\leq\,  \left(\,S_{F_{1}} \,\otimes\, S_{F_{2}}\,\right) \,\leq\, B\,D\,\left(\,I_{F} \,\otimes\, I_{G}\,\right)\\
&\Rightarrow\,A\,C\,I_{F \,\otimes\, G} \,\leq\, S_{F_{1} \,\otimes\, F_{2}} \,\leq\, B\,D\,I_{F \,\otimes\, G}. 
\end{align*}
This completes the proof.
\end{proof}

To each continuous frame in $n$-Hilbert space one can associate a dual continuous frame which is introduced as follows.\\

If \,$F_{1}$\, is a continuous frame associated to \,$\left(\,a_{\,2},\, \cdots,\, a_{\,n}\,\right)$\, for \,$H_{1}$\, with respect to \,$\left(\,X_{1},\, \mu_{1}\,\right)$\, and \,$F_{2}$\, is a continuous frame associated to \,$\left(\,b_{\,2},\, \cdots,\, b_{\,n}\,\right)$\, for \,$H_{2}$\, with respect to \,$\left(\,X_{2},\, \mu_{2}\,\right)$, then we may consider the dual continuous frame \,$G_{1}$\, associated to \,$(\,a_{\,2},\, \cdots,\, a_{\,n}\,)$\, of \,$F_{1}$\, and dual continuous frame \,$G_{2}$\, associated to \,$\left(\,b_{\,2},\, \cdots,\, b_{\,n}\,\right)$\, of \,$F_{2}$\, which satisfies the following: 
\begin{align}
&\left<\,f,\,g \,|\, a_{\,2},\, \cdots,\, a_{\,n}\,\right>_{\,1}\nonumber\\
& \,=\, \int\limits_{\,X_{1}}\,\left<\,f,\, F_{1}\,(\,x_{\,1}\,) \,|\, a_{\,2},\, \cdots,\, a_{\,n}\,\right>_{1}\,\left<\,G_{1}\,(\,x_{\,1}\,),\, g \,|\, a_{\,2},\, \cdots,\, a_{\,n}\,\right>_{1}\,d\mu_{\,1}\,(\,x_{1}\,),\label{eqn4.21}\\
&\left<\,f_{1},\,g_{1} \,|\, b_{\,2},\, \cdots,\, b_{\,n}\,\right>_{\,1}\nonumber\\
& \,=\, \int\limits_{\,X_{2}}\,\left<\,f_{1},\, F_{2}\,(\,x_{\,2}\,) \,|\, b_{\,2},\, \cdots,\, b_{\,n}\,\right>_{2}\,\left<\,G_{2}\,(\,x_{\,2}\,),\, g_{1} \,|\, b_{\,2},\, \cdots,\, b_{\,n}\,\right>_{2}\,d\mu_{\,2}\,(\,x_{2}\,),\label{eqn4.22} 
\end{align} 
for all \,$f,\, g \,\in\, H_{1}$\, and \,$f_{1},\, g_{1} \,\in\, H_{2}$.\\

Now, we give the definition of dual continuous frame in \,$H_{1} \,\otimes\, H_{2}$.
 
\begin{definition}
Let \,$\mathcal{F}$\, be a continuous frame associated to \,$(\,a_{\,2} \,\otimes\, b_{\,2},\, \cdots,\, a_{\,n} \,\otimes\, b_{\,n}\,)$\, for \,$H_{1} \,\otimes\, H_{2}$\, with respect to \,$(\,X,\, \mu\,)$.\,Then a frame \,$\mathcal{G}$\, associated to \,$(\,a_{\,2} \,\otimes\, b_{\,2},\, \cdots,\, a_{\,n} \,\otimes\, b_{\,n}\,)$\, satisfying 
\[f \,\otimes\, g \,=\, \int\limits_{\,X}\,\left<\,f \,\otimes\, g,\, \mathcal{F}\,(\,x\,) \,|\, a_{\,2} \,\otimes\, b_{\,2},\, \cdots,\, a_{\,n} \,\otimes\, b_{\,n}\,\right>\,\mathcal{G}\,(\,x\,)\,d\mu\,(\,x\,),\]for all \,$f \,\otimes\, g \,\in\, H_{1} \,\otimes\, H_{2}$, is called a dual continuous frame associated to \,$(\,a_{\,2} \,\otimes\, b_{\,2},\, \cdots,\, a_{\,n} \,\otimes\, b_{\,n}\,)$\, of \,$\mathcal{F}$.\,The pair \,$\left(\,\mathcal{F},\, \mathcal{G}\,\right)$\, is called a dual pair of continuous frames associated to \,$(\,a_{\,2} \,\otimes\, b_{\,2},\, \cdots,\, a_{\,n} \,\otimes\, b_{\,n}\,)$.    
\end{definition}

Next, we give a sufficient condition for two tensor product of continuous frames to form a pair of dual continuous frames in \,$H_{1} \,\otimes\, H_{2}$. 

\begin{theorem}\label{thm4.21}
Let \,$F_{1}$\, be a continuous frame associated to \,$\left(\,a_{\,2},\, \cdots,\, a_{\,n}\,\right)$\, for \,$H_{1}$\, with respect to \,$\left(\,X_{1},\, \mu_{1}\,\right)$\, and \,$F_{2}$\, is a continuous frame associated to \,$\left(\,b_{\,2},\, \cdots,\, b_{\,n}\,\right)$\, for \,$H_{2}$\, with respect to \,$\left(\,X_{2},\, \mu_{2}\,\right)$.\,Suppose \,$G_{1}$\, be the dual continuous frame associated to \,$\left(\,a_{\,2},\, \cdots,\, a_{\,n}\,\right)$\, of \,$F_{1}$\, and \,$G_{2}$\, be the dual continuous frame associated to \,$\left(\,b_{\,2},\, \cdots,\, b_{\,n}\,\right)$\, of \,$F_{2}$.\,Then \,$\mathcal{G} \,=\, G_{1} \,\otimes\, G_{2} \,:\, X \,\to\, H_{2} \,\otimes\, H_{2}$\, is a dual continuous frame associated to \,$\left(\,a_{\,2} \,\otimes\, b_{\,2},\, \cdots,\, a_{\,n} \,\otimes\, b_{\,n}\,\right)$\, for \,$H_{1} \,\otimes\, H_{2}$\, with respect to \,$(\,X,\, \mu\,)$\, of \,$\mathcal{F} \,=\, F_{1} \,\otimes\, F_{2} \,:\, X \,\to\, H_{1} \,\otimes\, H_{2}$.
\end{theorem}

\begin{proof}
By theorem \ref{th4.11}, \,$\mathcal{F} \,=\, F_{1} \,\otimes\, F_{2} \,:\, X \,\to\, H_{1} \,\otimes\, H_{2}$\, and \,$\mathcal{G} \,=\, G_{1} \,\otimes\, G_{2} \,:\, X \,\to\, H_{2} \,\otimes\, H_{2}$\, are continuous frames associated to \,$\left(\,a_{\,2} \,\otimes\, b_{\,2},\, \cdots,\, a_{\,n} \,\otimes\, b_{\,n}\,\right)$\, for \,$H_{1} \,\otimes\, H_{2}$\, with respect to \,$(\,X,\, \mu\,)$.\,Since \,$G_{1}$\, is a dual continuous frame associated to \,$\left(\,a_{\,2},\, \cdots,\, a_{\,n}\,\right)$\, of \,$F_{1}$\, and \,$G_{2}$\, is a dual continuous frame associated to \,$\left(\,b_{\,2},\, \cdots,\, b_{\,n}\,\right)$\, of \,$F_{2}$, for \,$f \,\in\, H_{1}$\, and \,$g \,\in\, H_{2}$, we have
\begin{align*}
& f\,=\, \int\limits_{\,X_{1}}\,\left<\,f,\, F_{1}\,(\,x_{\,1}\,) \,|\, a_{\,2},\, \cdots,\, a_{\,n}\,\right>_{1}\,G_{1}\,(\,x_{\,1}\,)\,d\mu_{\,1}\,(\,x_{1}\,),\\
& g\,=\, \int\limits_{\,X_{2}}\,\left<\,g,\, F_{2}\,(\,x_{\,2}\,) \,|\, b_{\,2},\, \cdots,\, b_{\,n}\,\right>_{2}\,G_{2}\,(\,x_{\,2}\,)\,d\mu_{\,2}\,(\,x_{2}\,). 
\end{align*}
Now, for each \,$f \,\otimes\, g \,\in\, H_{1} \,\otimes\, H_{2}$, we have
\begin{align*}
&f \,\otimes\, g\\
& \,=\, \int\limits_{\,X_{1}}\,\int\limits_{\,X_{2}}\,\left<\,f,\, F_{1}\,(\,x_{\,1}\,) \,|\, a_{\,2},\, \cdots,\, a_{\,n}\,\right>_{1}\,\left<\,g,\, F_{2}\,(\,x_{\,2}\,) \,|\, b_{\,2},\, \cdots,\, b_{\,n}\,\right>_{2}\mathcal{G}\,(\,x\,)d\mu\,(\,x\,) 
\end{align*}
where \,$\mathcal{G}\,(\,x\,) \,=\, G_{1}\,(\,x_{\,1}\,) \,\otimes\, G_{2}\,(\,x_{\,2}\,)$.\,By Fubini's theorem, we can write
\[f \,\otimes\, g \,=\, \int\limits_{\,X}\,\left<\,f \,\otimes\, g,\, \mathcal{F}\,(\,x\,) \,|\, a_{\,2} \,\otimes\, b_{\,2},\, \cdots,\, a_{\,n} \,\otimes\, b_{\,n}\,\right>\,\mathcal{G}\,(\,x\,)\,d\mu\,(\,x\,).\]
This completes the proof.       
\end{proof} 

In the following theorem, we will see that dual pair of continuous Bessel families is a dual pair of continuous frames in \,$H_{1} \,\otimes\, H_{2}$.

\begin{theorem}
Let \,$F_{1}$, \,$G_{1}$\, be the dual pair of continuous Bessel families associated to \,$(\,a_{\,2},\, \cdots,\, a_{\,n}\,)$\, for \,$H_{1}$\, with respect to \,$\left(\,X_{1},\, \mu_{1}\,\right)$\, having bounds \,$B_{1},\, B_{2}$\, and \,$F_{2}$, \,$G_{2}$\, be the dual pair of continuous Bessel families associated to \,$\left(\,b_{\,2},\, \cdots,\, b_{\,n}\,\right)$\, for \,$H_{2}$\, with respect to \,$\left(\,X_{2},\, \mu_{2}\,\right)$\, having bounds \,$D_{1},\, D_{\,2}$.\,Then \,$\mathcal{G} \,=\, G_{1} \,\otimes\, G_{2} \,:\, X \,\to\, H_{2} \,\otimes\, H_{2}$\, is a dual continuous frame associated to \,$\left(\,a_{\,2} \,\otimes\, b_{\,2},\, \cdots,\, a_{\,n} \,\otimes\, b_{\,n}\,\right)$\, for \,$H_{1} \,\otimes\, H_{2}$\, with respect to \,$(\,X,\, \mu\,)$\, of \,$\mathcal{F} \,=\, F_{1} \,\otimes\, F_{2} \,:\, X \,\to\, H_{1} \,\otimes\, H_{2}$.  
\end{theorem}

\begin{proof}
First, we show that \,$\mathcal{F} \,=\, F_{1} \,\otimes\, F_{2},\, \mathcal{G} \,=\, G_{1} \,\otimes\, G_{2}  \,:\, X \,\to\, H_{1} \,\otimes\, H_{2}$\, are continuous frames associated to \,$\left(\,a_{\,2} \,\otimes\, b_{\,2},\, \cdots,\, a_{\,n} \,\otimes\, b_{\,n}\,\right)$\, for \,$H_{1} \,\otimes\, H_{2}$\, with respect to \,$(\,X,\, \mu\,)$.\,Now, for each \,$f \,\otimes\, g \,\in\, H_{1} \,\otimes\, H_{2}$, using (\ref{eqn4.21}) and (\ref{eqn4.22}), we have
\begin{align*}
&\left\|\,f \,\otimes\, g,\, a_{\,2} \,\otimes\, b_{\,2},\, \cdots,\, a_{\,n} \,\otimes\, b_{\,n}\,\right\|^{\,2}\\
&=\,\left<\,f,\,f \,|\, a_{\,2},\, \cdots,\, a_{\,n}\,\right>_{\,1}\,\left<\,g,\,g \,|\, a_{\,2},\, \cdots,\, a_{\,n}\,\right>_{\,2}\\
&=\,\int\limits_{\,X_{1}}\,\left<\,f,\, F_{1}\,(\,x_{\,1}\,) \,|\, a_{\,2},\, \cdots,\, a_{\,n}\,\right>_{1}\,\left<\,G_{1}\,(\,x_{\,1}\,),\, f \,|\, a_{\,2},\, \cdots,\, a_{\,n}\,\right>_{1}\,d\mu_{\,1}\,(\,x_{1}\,)\,\times\\
& \hspace{1cm}\int\limits_{\,X_{2}}\,\left<\,g,\, F_{2}\,(\,x_{\,2}\,) \,|\, b_{\,2},\, \cdots,\, b_{\,n}\,\right>_{2}\,\left<\,G_{2}\,(\,x_{\,2}\,),\, g \,|\, b_{\,2},\, \cdots,\, b_{\,n}\,\right>_{2}\,d\mu_{\,2}\,(\,x_{2}\,)\\
&\leq\,\left(\,\int\limits_{\,X_{1}}\,\left|\,\left<\,f,\, F_{1}\,(\,x_{\,1}\,) \,|\, a_{\,2},\, \cdots,\, a_{\,n}\,\right>_{1}\,\right|^{\,2}d\mu_{\,1}\,(\,x_{1}\,)\,\right)^{1 \,/\, 2}\,\times\\
&\hspace{1cm}\left(\,\int\limits_{\,X_{1}}\,\left|\,\left<\,g,\, G_{1}\,(\,x_{\,1}\,) \,|\, a_{\,2},\, \cdots,\, a_{\,n}\,\right>_{1}\,\right|^{\,2}d\mu_{\,1}\,(\,x_{1}\,)\,\right)^{1 \,/\, 2}\,\times\\
&\hspace{1cm}\left(\,\int\limits_{\,X_{2}}\,\left|\,\left<\,g,\, F_{2}\,(\,x_{\,2}\,) \,|\, b_{\,2},\, \cdots,\, b_{\,n}\,\right>_{2}\,\right|^{\,2}d\mu_{\,2}\,(\,x_{2}\,)\,\right)^{1 \,/\, 2}\,\times\\
&\hspace{1cm}\left(\,\int\limits_{\,X_{2}}\,\left|\,\left<\,g,\, G_{2}\,(\,x_{\,2}\,) \,|\, b_{\,2},\, \cdots,\, b_{\,n}\,\right>_{2}\,\right|^{\,2}d\mu_{\,2}\,(\,x_{2}\,)\,\right)^{1 \,/\, 2}\\
&\leq\,\sqrt{B_{\,2}\,D_{\,2}}\,\left\|\,f,\, a_{\,2},\, \cdots,\, a_{\,n}\,\right\|_{1}\,\left\|\,g,\, b_{\,2},\, \cdots,\, b_{\,n}\,\right\|_{2}\,\times\\
&\hspace{1cm}\left(\int\limits_{\,X}\,\left|\,\left<\,f \,\otimes\, g,\, \mathcal{F}\,(\,x\,) \,|\, a_{\,2} \,\otimes\, b_{\,2},\, \cdots,\, a_{\,n} \,\otimes\, b_{\,n}\,\right>\,\right|^{\,2}\,d\mu\,(\,x\,)\right)^{1 \,/\, 2}\\
&\Rightarrow\,\dfrac{1}{B_{\,2}\,D_{\,2}}\,\left\|\,f \,\otimes\, g,\, a_{\,2} \,\otimes\, b_{\,2},\, \cdots,\, a_{\,n} \,\otimes\, b_{\,n}\,\right\|^{\,2}\\
&\hspace{1cm}\leq\,\int\limits_{\,X}\,\left|\,\left<\,f \,\otimes\, g,\, \mathcal{F}\,(\,x\,) \,|\, a_{\,2} \,\otimes\, b_{\,2},\, \cdots,\, a_{\,n} \,\otimes\, b_{\,n}\,\right>\,\right|^{\,2}\,d\mu\,(\,x\,).
\end{align*} 
Thus, \,$\mathcal{F}$\, is a continuous frame associated to \,$\left(\,a_{\,2} \,\otimes\, b_{\,2},\, \cdots,\, a_{\,n} \,\otimes\, b_{\,n}\,\right)$\, for \,$H_{1} \,\otimes\, H_{2}$\, with respect to \,$(\,X,\, \mu\,)$\, having bounds \,$\dfrac{1}{B_{\,2}\,D_{\,2}}$\, and \,$B_{\,1}\,D_{\,1}$.\,Similarly, it can be shown that \,$\mathcal{G}$\, is a continuous frame associated to \,$\left(\,a_{\,2} \,\otimes\, b_{\,2},\, \cdots,\, a_{\,n} \,\otimes\, b_{\,n}\,\right)$\, for \,$H_{1} \,\otimes\, H_{2}$\, with respect to \,$(\,X,\, \mu\,)$.\,Now, by theorem \ref{thm4.21}, \,$\left(\,\mathcal{F},\, \mathcal{G}\,\right)$\, is a dual pair of continuous frames associated to \,$(\,a_{\,2} \,\otimes\, b_{\,2},\, \cdots,\, a_{\,n} \,\otimes\, b_{\,n}\,)$.\,This completes the proof.    
\end{proof} 

Now, we end this section by discussing the idea of continuous Bessel multiplier in \,$H_{1} \,\otimes\, H_{2}$.   

\begin{definition}
Let \,$\mathcal{F}$\, and \,$\mathcal{G}$\, be continuous Bessel families associated to \,$(\,a_{\,2} \,\otimes\, b_{\,2},\, \cdots,\, a_{\,n} \,\otimes\, b_{\,n}\,)$\, for \,$H_{1} \,\otimes\, H_{2}$\, with respect to \,$\left(\,X,\, \mu\,\right)$\, having bounds \,$B_{1}$\, and \,$B_{2}$\, and \,$m \,:\, X \,\to\, \mathbb{C}$\, be a measurable function.\,The operator \,$\mathcal{M}_{m,\, \mathcal{F},\, \mathcal{G}} \,:\, H_{F} \,\otimes\, H_{G} \,\to\, H_{F} \,\otimes\, H_{G}$\, defined by
\begin{align}
&\mathcal{M}_{m,\, \mathcal{F},\, \mathcal{G}}\,(\,f \,\otimes\, g\,)\nonumber\\
& \,=\,: \int\limits_{\,X}\,m\,(\,x\,)\,\left<\,f \,\otimes\, g,\, \mathcal{F}\,(\,x\,) \,|\, a_{\,2} \,\otimes\, b_{\,2},\, \cdots,\, a_{\,n} \,\otimes\, b_{\,n}\,\right>\,\mathcal{G}\,(\,x\,)\,d\mu\,(\,x\,),\label{eqn4.23}
\end{align} 
for all \,$f \,\otimes\, g \,\in\, H_{F} \,\otimes\, H_{G}$, is called continuous Bessel multiplier associated to \,$(\,a_{\,2} \,\otimes\, b_{\,2},\, \cdots,\, a_{\,n} \,\otimes\, b_{\,n}\,)$\, of \,$\mathcal{F}$\, and \,$\mathcal{G}$\, with respect to \,$m$. 
\end{definition}

\begin{note}
Let \,$F_{1}$, \,$G_{1}$\, be continuous Bessel families associated to \,$(\,a_{\,2},\, \cdots,\, a_{\,n}\,)$\, for \,$H_{1}$\, with respect to \,$\left(\,X_{1},\, \mu_{1}\,\right)$\, and \,$F_{2}$, \,$G_{2}$\, be continuous Bessel families associated to \,$(\,b_{\,2},\, \cdots,\, b_{\,n}\,)$\, for \,$H_{2}$\, with respect to \,$\left(\,X_{2},\, \mu_{2}\,\right)$\, and \,$m_{1} \,:\, X_{1} \,\to\, \mathbb{C}$, \,$m_{2} \,:\, X_{2} \,\to\, \mathbb{C}$\, be two measurable function.\,Suppose \,$M_{m_{1},\, F_{1},\, G_{1}} \,:\, H_{F} \,\to\, H_{F}$\, be a continuous Bessel multiplier associated to \,$(\,a_{\,2},\, \cdots,\, a_{\,n}\,)$\, of \,$F_{1}$\, and \,$G_{1}$\, with respect to \,$m_{1}$\, and \,$M_{m_{2},\, F_{2},\, G_{2}} \,:\, H_{G} \,\to\, H_{G}$\, be a continuous Bessel multiplier associated to \,$(\,b_{\,2},\, \cdots,\, b_{\,n}\,)$\, of \,$F_{2}$\, and \,$G_{2}$\, with respect to \,$m_{2}$.\,Now, by theorem \ref{th4.11}, \,$\mathcal{F} \,=\, F_{1} \,\otimes\, F_{2},\, \mathcal{G} \,=\, G_{1} \,\otimes\, G_{2}  \,:\, X \,\to\, H_{1} \,\otimes\, H_{2}$\, are continuous Bessel families associated to \,$\left(\,a_{\,2} \,\otimes\, b_{\,2},\, \cdots,\, a_{\,n} \,\otimes\, b_{\,n}\,\right)$\, for \,$H_{1} \,\otimes\, H_{2}$\, with respect to \,$(\,X,\, \mu\,)$.\,From (\ref{eqn4.23}), for each \,$f \,\otimes\, g \,\in\, H_{F} \,\otimes\, H_{G}$, we can write
\begin{align*}
&\mathcal{M}_{m,\, \mathcal{F},\, \mathcal{G}}\,(\,f \,\otimes\, g\,)\\
& \,=\,: \int\limits_{\,X}\,m\,(\,x\,)\,\left<\,f \,\otimes\, g,\, \mathcal{F}\,(\,x\,) \,|\, a_{\,2} \,\otimes\, b_{\,2},\, \cdots,\, a_{\,n} \,\otimes\, b_{\,n}\,\right>\,\mathcal{G}\,(\,x\,)\,d\mu\,(\,x\,)\\
&\,=\,: \int\limits_{\,X_{1}}\,m_{1}\,(\,x_{1}\,)\,\left<\,f,\, F_{1}\,(\,x_{\,1}\,) \,|\, a_{\,2},\, \cdots,\, a_{\,n}\,\right>_{1}\,G_{1}\,(\,x_{\,1}\,)\,d\mu_{\,1}\,(\,x_{1}\,)\,\,\otimes\\
&\hspace{1cm}\int\limits_{\,X_{2}}\,m_{2}\,(\,x_{2}\,)\,\left<\,g,\, F_{2}\,(\,x_{\,2}\,) \,|\, b_{\,2},\, \cdots,\, b_{\,n}\,\right>_{2}\,G_{2}\,(\,x_{\,2}\,)\,d\mu_{\,2}\,(\,x_{2}\,)\\
&=\,:\, M_{m_{1},\, F_{1},\, G_{1}}\,f \,\otimes\, M_{m_{2},\, F_{2},\, G_{2}}\,g \,=\, \left(\, M_{m_{1},\, F_{1},\, G_{1}} \,\otimes\, M_{m_{2},\, F_{2},\, G_{2}}\,\right)\,(\,f \,\otimes\, g\,).
\end{align*}
Thus, \,$\mathcal{M}_{m,\, \mathcal{F},\, \mathcal{G}} \,=\, M_{m_{1},\, F_{1},\, G_{1}} \,\otimes\, M_{m_{2},\, F_{2},\, G_{2}}$. 
\end{note}

\begin{remark}
According to the theorem \ref{thm4.22}, the continuous Bessel multiplier associated to \,$(\,a_{\,2} \,\otimes\, b_{\,2},\, \cdots,\, a_{\,n} \,\otimes\, b_{\,n}\,)$\, of \,$\mathcal{F}$\, and \,$\mathcal{G}$\, with respect to \,$m$\, is well defined and bounded. 
\end{remark}

\end{document}